\documentclass[11pt,a4paper,twoside,reqno]{amsart}
\usepackage{amsmath}
\usepackage{amsfonts}
\usepackage{amssymb}
\usepackage{amstext}
\usepackage{amsbsy}

\usepackage{graphicx}

\usepackage{epsf}
\newcommand{\Proof}{\noindent \textbf{Proof.} }

\theoremstyle{plain}
\begingroup

\newtheorem{thm}{Theorem}
\newtheorem{prop}[thm]{Proposition}

\newtheorem{lem}[thm]{Lemma}

\endgroup
\newtheorem{rem}[thm]{\textbf{Remark}}

\theoremstyle{definition}

\newcommand{\ch}{\cosh}
\newcommand{\sh}{\sinh}

\newcommand{\tg}{\tan}

\newcommand{\pt}{\frac{\partial}{ \partial t}}
\newcommand{\px}{\frac{\partial}{ \partial x}}
\newcommand{\py}{\frac{\partial}{ \partial y}}
\newcommand{\pz}{\frac{\partial}{ \partial z}}

\newcommand{\ov}[1]{\overline{#1}}
\newcommand{\wt}[1]{\widetilde{#1}}

\newcommand{\hi}[1]{\mathbb{H}^#1}
\newcommand{\m}[1]{\mathbb{R}^#1}

\newcommand{\Z}{\mathbb{Z}}
\newcommand{\Sp}{\mathbb{S}}
\newcommand{\M}{\mathbb{M}}
\newcommand{\Mt}{\mathbb{M}^3}
\newcommand{\cm}{\mathcal{C}}

\newcommand{\R}{\mathbb{R}}
\newcommand{\sd}{\mathbb{S}^2}

\newcommand{\hd}{\mathbb{H}^2}
\newcommand{\h}{_\mathbb{H}}

\newcommand{\C}{\mathbb{C}}

\def\rmd{\mathop{\rm d\kern -1pt}\nolimits}

\begin{document}
\centerline{\large\bf Totally umbilic surfaces in homogeneous 3-manifolds}

\vskip4mm
\centerline{\textbf{\textit{Rabah Souam}} $\&$
\textbf{\textit{Eric Toubiana}}}
\title{ }

\address{ Institut de Math{\'e}matiques de Jussieu\newline
 \indent CNRS UMR 7586 - Universit{\'e} Paris Diderot - Paris 7\newline
 \indent G\'eom\'etrie et Dynamique \newline
 \indent Site Chevaleret \newline
\indent Case 7012       \newline  
 \indent         75205 - Paris Cedex 13, France \newline
\indent  email: souam@math.jussieu.fr\newline
\indent  email: toubiana@math.jussieu.fr}



\begin{abstract}
We discuss existence and classification of totally umbilic surfaces in
the model geometries of Thurston and the Berger spheres. We classify 
such surfaces in $\hd\times\R$, $\sd\times \R$ and the $Sol$ group. We
prove nonexistence in the Berger spheres and in the remaining
model geometries other than the space forms.
\end{abstract}


\date{\today}

\maketitle

\vskip8mm

\section{Introduction}\label{Sec.Intro}

During the recent years, there has been a rapidly growing interest in
the geometry of surfaces in $\sd \times \R$ and $ \hd \times \R$
focusing on minimal and constant mean curvature surfaces. This was
initiated by H. Rosenberg, \cite{[Rosenberg]}. More generally many works
are devoted to studying the geometry of surfaces in homogeneous 
3-manifolds. See for example 
\cite{[M-R]}, \cite{[Daniel1]},
 \cite{[Daniel]}, \cite{[Sanini]}, 
\cite{[Nelli-R]}, \cite{[Hauswirth]}, 
\cite{[Hauswirth1]}, 
\cite{[Pedrosa]}, \cite{[Earp-T]}, 
\cite{[Galvez]}, \cite{[Abresch-R]},
 \cite{[Isabel]}, 
\cite{[Rato-CMC]} and \cite{[Earp]}. 

In the space forms the classification of totally umbilic surfaces 
is well known and very useful, see \cite{[Spivak]}. In $\R^3$ they are 
planes and round
spheres and in $\mathbb{S}^3$ they are round spheres. In $\hi3$ they
are totally geodesic planes and their equidistants, horospheres
and round spheres. In particular they all have constant mean curvature.

A natural question is to understand the totally
umbilic surfaces in the remaining homogeneous 3-manifolds. Untill now
the only known result in this direction was the non-existence of totally
umbilic surfaces in the Heisenberg space due to A. Sanini, see
\cite{[Sanini]}. In this paper we study totally umbilic surfaces in
simply connected and homogeneous 3-manifolds. More precisely we first 
consider the manifolds having a 4-dimensional 
isometry group, denoted by 
$\Mt (\kappa,\tau)$ (see section \ref{nonumbilic}).
Namely these manifolds are 
$\sd (\kappa) \times \R$ ($\kappa>0,\ \tau=0$), 
$\hd (\kappa)\times \R$ ($\kappa<0,\ \tau=0$), 
the Berger spheres \newline
 ($\kappa>0,\ \tau\not=0$)
and the manifolds having the isometry group of either
the Heisenberg space ($\kappa=0,\ \tau\not=0$)
 or $\wt{\mathrm{PSL}_2(\R)}$ ($\kappa<0,\ \tau\not=0$), see
 \cite{[Bonahon]}, \cite{[Scott]} or
\cite{[Thurston]}. Except for the Berger spheres, these manifolds are
four of the eight model geometries of Thurston \cite{[Thurston]}. The
remaining model geometries are the three space forms and the $Sol$
geometry which has a 3-dimensional isometry group. As a matter of fact
we also consider the $Sol$ geometry. 

 In section  \ref{nonumbilic} we prove (Theorem \ref{nonexistence})
the non-existence of totally umbilic surfaces, in particular the
non-existence of totally geodesic ones, in the homogeneous 
manifolds $\Mt (\kappa,\tau)$ for $\tau\not=0$, that is those 
which are not Riemannian products. This extends the result of
 Sanini, \cite{[Sanini]}. 

 In section \ref{sphere} we construct and 
classify all rotational and totally umbilic surfaces in 
$\sd (\kappa)\times \R$. In section \ref{plan hyperbolique} we construct and 
classify all totally umbilic surfaces in $\hd (\kappa)\times \R$ 
which are invariant under 
a one-parameter group of ambient isometries. Except for the totally
geodesic ones, these surfaces do not have constant mean curvature.
In section 
\ref{unicite}, we prove that the surfaces obtained in sections
\ref{sphere} and \ref{plan hyperbolique} are the only 
totally umbilic surfaces in respectively $\sd(\kappa) \times \R$ and 
$\hd (\kappa)\times \R$. 

In section \ref{Sol} we show that there exist,
up to ambient isometries, only two
totally umbilic surfaces in $Sol$, one of them being totally geodesic.

Finally, in section \ref{conform} we apply our results to prove that any conformal diffeomorphism of 
$\hd\times\R,$ $\sd\times\R$ and {\it Sol} is an isometry.

 Throughout this paper all the surfaces are assumed of class
 $C^3$, see however the remark \ref{regularite}.

 We are grateful to H. Rosenberg for valuable comments and to the referee for
his observations which improved the paper.

\section{Non-existence of totally umbilic surfaces in some 
homogeneous 3-manifolds}\label{nonumbilic}
In this section we consider the connected and simply connected 
homogeneous Riemannian 3-manifolds, whose isometry group has dimension 4 
and which are not Riemannian products.
 We recall that such a manifold is a fibration over some
complete and simply connected
surface, $\M^2 (\kappa)$, of constant curvature $\kappa \in \R$, 
with geodesic fibers.
Actually, for each $\kappa$, there is a one-parameter family
$\Mt (\kappa,\tau)$ of such fibrations, parametrized by 
the bundle curvature $\tau \in \R ^\ast$. 
The unit vector field $\xi$ tangent to the fibers is a 
Killing field and satisfies :
\begin{equation}\label{vectoriel}
\ov \nabla_X \xi =\tau (X\wedge \xi),
\end{equation}
for any tangent vector $X$ in $T\Mt (\kappa,\tau)$, where
$\ov \nabla$ is the connection on $\Mt (\kappa,\tau)$. 
The field $\xi$ defines the \textsl{vertical} direction of the
Riemannian submersion \newline
$\Mt (\kappa,\tau) \rightarrow \M^2 (\kappa)$. 
As a matter of fact, the bundle curvature $\tau$ can be zero, but in this
case $\Mt (\kappa,0)$ is just a Riemannian product 
$\M^2(\kappa) \times \R$. These product manifolds will be considered in the
following sections. 
Moreover we assume $\kappa -4\tau^2\not=0$, otherwise the 
manifold is a space form and its isometry group has dimension 6.
These manifolds are of three types :
when $\kappa >0$ 
they are the Berger spheres, for
$\kappa=0$ they have the isometry group of the 
Heisenberg space, Nil$_3$, and  
for $\kappa <0$ they have the isometry group of 
$\wt{\mathrm{PSL}_2(\R)}$.

For more details we refer to \cite{[Bonahon]}, 
\cite{[Scott]} and \cite{[Thurston]}.

We can now state our result.

\begin{thm}\label{nonexistence}
There exist no totally umbilic surfaces (even non complete) in 
the 3-manifolds $\Mt (\kappa,\tau)$, with $\tau\not= 0$ and 
$\kappa -4\tau^2\not=0$. In particular, there are no totally 
geodesic surfaces.
\end{thm} 

For the special case of the Heisenberg space (
$\kappa = 0, \tau=1/2$), we 
recover the result proved by A. Sanini, 
see \cite{[Sanini]}.

\Proof
 Let $S$ be an immersed totally umbilic surface in 
 $\Mt (\kappa,\tau)$. 
Locally $S$ is the image of an embedding 
$X:\Omega \rightarrow \Mt (\kappa,\tau)$, where $\Omega$ is an open disk
in $\m 2$. Call $(u,v)$ the coordinates on $\Omega$ and consider
a unit normal field $N$ on $X(\Omega)$. As $X$ is
totally umbilic, there exists a function 
$\lambda :\Omega \rightarrow \R$ such that 
\begin{equation*}
\left\{ 
\begin{aligned}
\ov{\nabla}_{X_u}N =\lambda X_u \\
\ov{\nabla}_{X_v}N =\lambda X_v
\end{aligned} \right.
\end{equation*}
Therefore
\begin{equation*}
\left\{ 
\begin{aligned}
\ov{\nabla}_{X_v}(\ov{\nabla}_{X_u}N) = 
\lambda_v X_u +\lambda \ov{\nabla}_{X_v}X_u\\
\ov{\nabla}_{X_u}(\ov{\nabla}_{X_v}N) = 
\lambda_u X_v +\lambda \ov{\nabla}_{X_u}X_v
\end{aligned} \right.
\end{equation*}

Substracting the second equation from the first one we get
\begin{equation*}
\ov{\nabla}_{X_v}(\ov{\nabla}_{X_u}N)-
\ov{\nabla}_{X_u}(\ov{\nabla}_{X_v}N)=
\lambda_v X_u - \lambda_u X_v.
\end{equation*}

That is 
\begin{equation}\label{tenseur}
R(X_u,X_v)N=\lambda_v X_u - \lambda_u X_v,
\end{equation}
where $R$ denotes the curvature tensor of $\Mt (\kappa,\tau)$.

We define the function $\nu$
on $\Omega$ setting $\nu := \langle N,\xi \rangle$. We denote by $T$ the
projection of $\xi$ on $S$, that is $T=\xi-\nu N$.

As the projection $\Mt (\kappa,\tau)\rightarrow \M^2 (\kappa)$ is a 
Riemannian submersion, we have the following formula derived by 
Daniel, see \cite{[Daniel]}:
\begin{equation*}
R(X_u,X_v)N=(\kappa -4\tau^2)\nu 
\left( \langle X_v,T\rangle X_u - \langle X_u,T\rangle X_v\right).
\end{equation*}
Taking into account the relation (\ref{tenseur}) we get
\begin{equation}\label{gradient}
\nabla \lambda =(\kappa -4\tau^2)\nu T,
\end{equation}
where $\nabla$ denotes the gradient on $S$.

Observe that if $T=0$ on a nonempty open set, then we can take 
$N=\xi$ on this set and we deduce from (\ref{vectoriel}) that this
surface cannot be umbilic. We
can thus assume that $T$ does not vanish on $\Omega$.

Set $JT=N\wedge T$, thus $JT$ is tangent to $S$ and horizontal.

\medskip 
\textbf{Claim:} We have $[T,JT]\equiv 0$.

We need to show that:  $\ov \nabla_T JT=\ov \nabla_{JT} T$.
Since $JT=N\wedge T=N\wedge (\xi-\nu N)=N\wedge \xi$, we have 
using (\ref{vectoriel})
\begin{equation*}
\ov \nabla_T JT=\ov \nabla_T (N\wedge \xi)=\ov \nabla_TN \wedge \xi 
+N\wedge \ov \nabla_T\xi=\lambda(T\wedge \xi)+\tau (N\wedge(T\wedge \xi)).
\end{equation*}
As $T\wedge \xi= T\wedge(T+\nu N)=\nu(T\wedge N)=-\nu JT$, we deduce
that 
\begin{equation*}
\ov \nabla_T JT=-\lambda \nu JT +\tau \nu T.
\end{equation*}
On the other hand
\begin{equation}\label{Eq1}
\ov \nabla _{JT} T=\ov \nabla _{JT}(\xi -\nu N)=\tau (JT\wedge \xi)-
JT.(\nu) N-\lambda \nu JT.
\end{equation}  
We have
\begin{equation}\label{Eq2}
JT \wedge \xi =JT \wedge (T+\nu N)=JT \wedge T + \nu JT \wedge N=
-\vert T \vert^2 N + \nu T,
\end{equation}
and 
\begin{equation*}
\begin{aligned}
JT.(\nu) =&\ JT. \langle N,\xi\rangle=
\langle \ov \nabla _{JT}N,\xi\rangle  +
\langle N,\ov \nabla _{JT}\xi\rangle \\
=&\ \lambda\langle JT, \xi\rangle +\tau \langle N,JT\wedge \xi\rangle
=\tau \langle N,JT\wedge \xi\rangle.
\end{aligned}
\end{equation*}
Using (\ref{Eq2}) we obtain 
\begin{equation}\label{Eq3}
JT.(\nu)=-\tau \vert T \vert^2.
\end{equation}
Inserting (\ref{Eq2}) and (\ref{Eq3}) in (\ref{Eq1}) we end with:
\begin{equation*}
\ov \nabla _{JT} T=\tau \nu T -\lambda \nu JT =
\ov \nabla _T JT,
\end{equation*}
which proves the claim.


Now (\ref{gradient}) implies
\begin{equation*}
\left\{ \begin{aligned}
JT.(\lambda)= &\ 0  \\
T.(\lambda)=&\ (\kappa-4\tau^2)\nu \vert T\vert^2= (\kappa-4\tau^2)\nu 
(1-\nu^2)\\
\end{aligned}\right.
\end{equation*}
Since $[T,JT]=0$, we get 
$(\kappa-4\tau^2)JT.(\nu-\nu^3)=0$. As $\kappa-4\tau^2\not= 0$ we
infer that $(1-3\nu^2)JT.(\nu) =0$. This implies easily 
$$
JT.(\nu)=0.
$$
As $JT.(\nu)=-\tau \vert T\vert^2$, see (\ref{Eq3}), and $\tau\not= 0$, 
we deduce that $T\equiv 0$, which is a contradiction. This concludes
the proof.\qed

\section{Symmetric totally umbilic surfaces in 
$\sd (\kappa) \times \R$}\label{sphere}
In this section we classify totally umbilic surfaces which are
rotationally invariant in $\sd \times \R$. The classification in
$\M ^2 (\kappa)\times \R$, for any $\kappa >0$ is completely analogous.
 We will see that besides the obvious totally
geodesic ones, up to isometries of 
$\sd \times \R$, 
there are two one-parameter families of complete 
totally umbilic rotationally invariant surfaces  
 homeomorphic to the two-sphere 
and a unique complete surface 
which has the topology of $\m 2$.
The surfaces of the
first family are homologous to zero and those of the second family are not.
Moreover these surfaces are embedded, analytic and any totally umbilic 
rotationally invariant
surface is a part of one of these complete surfaces.

A rotational surface in $\sd \times \R$ is by definition a surface
obtained by rotating a curve in a totally geodesic cylinder 
$\cm:=\Gamma \times \R$, where $\Gamma \subset \sd$ is a geodesic, around 
an axis $R=\{p\}\times \R$ where $p$ is a fixed point of $\Gamma$.

In the coordinates $(x,y,t)$ 
given by the stereographic 
projection with respect to the north pole, the metric on 
$\sd \times \R$ reads as follows:
$$
\rmd \tilde{s}^2=\left( \dfrac{2}{1+x^2+y^2}\right)^2
(\rmd x^2+\rmd y^2) + \rmd t^2,
$$
where $x,y,t \in \R$.

 Up to an ambient isometry we can assume
that $\Gamma \subset \sd$ corresponds to the complete 
geodesic defined by $y=0$ and that  $p=(0,0,0)\in \Gamma$ is 
the south pole of $\sd$. Therefore the axis
is $R=\{(0,0,t),\ t\in \R\}$.

\begin{rem}\label{poles}
Let us remark that for any given curve $\alpha$ in the cylinder $\cm$, the
surface generated by rotating $\alpha$ around the axis through the
south pole, that is $R$, is the same as the one generated    
by rotating $\alpha$ around the axis through the
north pole.
\end{rem}

We consider the
vertical  (noncomplete) geodesic plane $P=\{y=0\}\subset \cm$, that is 
$\cm =P\cup (\{N\}\times \R)$, where $N\in \sd$ is the north pole.

 Let $\rho \in \mathopen]-\pi,\pi\mathclose[$ 
denote the signed distance to the origin $(0,0)$ 
on $\Gamma \cap P$. Thus we have 
$x=\tg(\rho/2)$. In the coordinates 
$(\rho,t)$ the metric on the plane $P$ writes 
$$
\rmd s^2=\rmd \rho^2 + \rmd t^2.
$$

Consider now a smooth curve $\alpha (s)=(\rho (s),t(s))$
parametrized
by arclength in $P$. Let $\theta (s)$ be the oriented angle between
the $\rho$-axis and $\alpha^\prime (s)$.
Therefore, we have:
\begin{equation}\label{tangent}
\left\{ 
\begin{aligned}
\rho^\prime (s)=&\ \cos \theta(s)\\
t^\prime(s)=&\ \sin\theta(s)\\
\end{aligned}\right.
\end{equation}

In the plane $P$ we consider the unit normal $N$ to the curve 
$\alpha$ so that the basis ($\alpha^\prime (s),N(s))$ is positively oriented for
each $s$. We orient by $N$ the symmetric
surface generated by $\alpha$. The principal curvatures 
computed with respect to this orientation are as follows:

\begin{equation*}
\left\{ 
\begin{aligned}
\lambda_1(s) =&\  \theta ^\prime (s) \\
\lambda_2(s) =&\ \dfrac{\sin \theta (s)}{\tan \rho(s)}.\\
\end{aligned}\right.
\end{equation*}
Thus, the umbilicity condition is:

\begin{equation}\label{umbilic}
\theta^\prime (s)=\dfrac{\sin\theta(s)}{\tan \rho (s)}\cdot
\end{equation}

A priori the equation (\ref{umbilic}) makes sense only for 
$\rho \not=0$, but  as we will see later the surfaces we obtain 
are regular even at such points.  

Differentiating the first equation in (\ref{tangent}) and using
equation (\ref{umbilic}) we get

\begin{equation}\label{second}
\rho^{\prime \prime}=\dfrac{(\rho^{\prime 2}-1)}{\tan \rho}\cdot
\end{equation}

Assume that $\rho^\prime (s_0)=1$ for some $s_0$ where 
$\rho (s_0)\not=0$. Note that the
function $f(s)=\rho(s_0)+s-s_0$ is a solution of (\ref{second}) with
the same initial conditions at $s_0$ than $\rho$. Therefore 
$\rho \equiv f$ and $t^\prime \equiv 0$ and the surface is part of a
slice $\sd \times \{t_0\}$. The same happens in case where 
$\rho^\prime (s_0)=-1$. Henceforth we will assume that 
$\rho^{\prime 2}(s)\not= 1$ for all $s$ and (\ref{second}) is
equivalent to
$$
\dfrac{\rho^{\prime \prime}}{(\rho^{\prime 2}-1)}=
\dfrac{\cos \rho}{\sin \rho} \cdot
$$
Multiplying both sides by $2\rho^\prime$ and integrating we get
$$
\rho^{\prime 2}-1=\lambda\sin^2\rho,
$$
for some nonzero real constant $\lambda$. Since the curve $\alpha$ is
parametrized by arclength we must have $\rho^{\prime 2}<1$. 
Thus $\lambda=-a^2$ for some $a>0$.

Conversely, any solution $\rho $ of (\ref{second}) satisfying 
$\rho^{\prime 2}\leq1$ defines a function 
$\theta (s)$ setting $\cos \theta (s)=\rho^\prime (s)$.
Consider the function $t$ defined by setting 
$t^\prime (s)=a\sin \rho(s)$ and $t(s_0)=t_0$ for some  $s_0$ in the
domain of $\rho$ and some real number $t_0$. Then $t$ satisfies the
second equation of (\ref{tangent}) and therefore the curve 
$\alpha(s)=(\rho(s),t(s))\in \cm$ generates a rotational totally
umbilic surface in $\sd \times \R$.

Let $\rho$ be a solution of (\ref{second})
satisfying $\rho^{\prime 2}<1$. Observe that
equations (\ref{tangent}) and (\ref{umbilic}) show that $\rho^\prime$ cannot be 
identically zero on an open interval unless the generated surface is
part of the totally geodesic cylinder $\rho \equiv \pi/2$.
 Henceforth we
assume that $\rho$ is not this trivial solution and so, 
up to restricting the domain of $\rho$, we can suppose that $\rho$ takes its
values in $\mathopen]0,\pi/2\mathopen[$ or 
$\mathopen]\pi/2,\pi\mathopen[$.  

 So, we can consider an interval on
which $\rho^\prime$ never vanishes.
Changing $s$ into $-s$ if needed we can suppose that 
$ \rho^\prime >0$. 
Therefore we get
\begin{equation}\label{intprem}
\rho^\prime =\sqrt{1-a^2\sin^2\rho}.
\end{equation}

It is interesting to note that, when $a\not=1$, 
the function $\rho$ is
the {\it Jacobi amplitude} function : $\rho (s)=\mathrm{am}(s,a^2)$ and, up to
the
sign, we have $t^\prime (s)= a \ \mathrm{sn}(s, a^2)$ and 
$\theta^\prime (s)= a\  \mathrm{cn}(s,a^2)$, where 
$\mathrm{sn}(s,a^2)=\sin \mathrm{am}(s,a^2)$ and 
$\mathrm{cn}(s,a^2)=\cos \mathrm{am}(s,a^2)$ are respectively the {\it sinus and
cosinus amplitudinis} elliptic Jacobi functions, see for instance 
\cite[Chapter 16]{[Abra-Ste]} and \cite[pp 286-307]{[San-Ger]}.
However, for reader's convenience and to be self-contained, 
we will treat in a direct and elementary way this ODE.

Now observe that the transformation $(\rho,t)\mapsto (\pi -\rho,t)$ is
an isometry which changes rotations around the axis through the
south pole into rotations around the axis through the
north pole. Therefore, taking into account Remark \ref{poles},
we can assume that, up to an ambient isometry,
$\rho$ takes its
values in $\mathopen]0,\pi/2\mathopen[$.

Let us call $\rho_a$ the maximal solution of (\ref{intprem}) extending
$\rho$ without restrictions on its values, that is for the moment 
we do not require that $\rho_a$ takes its values in  
$\mathopen]-\pi,\pi\mathopen[$.

\begin{lem}\label{zero}
Up to a reparametrization of the form $s\rightarrow s+c$ for some real
constant $c$, the maximal solution $\rho_a$
is defined on an interval
$\mathopen]-\delta,\delta\mathopen[$, where 
$\delta \in \mathopen]0,+\infty\mathopen]$. Furthermore 
$\rho_a$ is odd and so satisfies $\rho_a (0)=0$ and 
$\rho_a ^\prime (0)=1$.
\end{lem}

\Proof Let us call 
$\mathopen]u,v\mathopen[$ the domain of $\rho_a$ where 
$-\infty \leq u<v \leq +\infty$. 

We first show that 
$\rho_a$ vanishes at some point. 
Since $\rho_a^\prime \leq 1$ such a point clearly exists if
$u=-\infty$. In case $u$ is finite, $\rho_a$ has a limit 
$l \in \mathopen[-\infty,+\infty\mathopen[$
at $u$ as it is nondecreasing. 

If $l<0$ then 
$\rho_a$ vanishes at some point since $\rho$ is positive.

Consider now the case where $l\geq 0$
and call
$I\subset   \mathopen]u,v\mathopen[$ the domain of $\rho$.
Suppose that, as $s$ decreases starting from $I$, $\rho_a$ never
vanishes. Then the function $\rho_a^\prime$ increases 
(since $\rho (I)\subset \mathopen]0,\pi/2\mathopen[$) 
and thus 
has a positive limit at $u$. But then we could extend the solution
$\rho_a$ of equation (\ref{intprem}) beyond $u$, which 
contradicts the maximality of $\rho_a$. 

Therefore $\rho_a (s_0)=0$ for some point $s_0$. Changing $s$ into
$s-s_0$ we can assume that $\rho_a (0)=0$. The function 
$f(s):=-\rho_a(-s)$ is then also a solution of (\ref{intprem}) satisfying
$f(0)=0$. We conclude that $f\equiv \rho_a$, which means $\rho_a$
is odd. \qed

\begin{lem}\label{a<1}
Suppose $a\in \mathopen]0,1\mathopen[$. Then $\rho_a$ is defined on the
whole of $\R$ and it gives rise to a unique, up to an ambient
isometry, curve $\alpha_a$ generating a rotational totally 
umbilic surface. This curve is an analytic 
Jordan curve in the
cylinder $\cm$, it is nonhomologous to zero and symmetric with
respect to the axis of rotation $R$. 
The rotational
totally umbilic surface, $S_1(a)$,
generated by $\alpha_a$ is analytic, embedded and 
homeomorphic to the sphere. 
Moreover $S_1(a)$ is nonhomologous to zero in $\sd \times \R$.
\end{lem}

\begin{figure}[ht]
\includegraphics[scale=0.45]{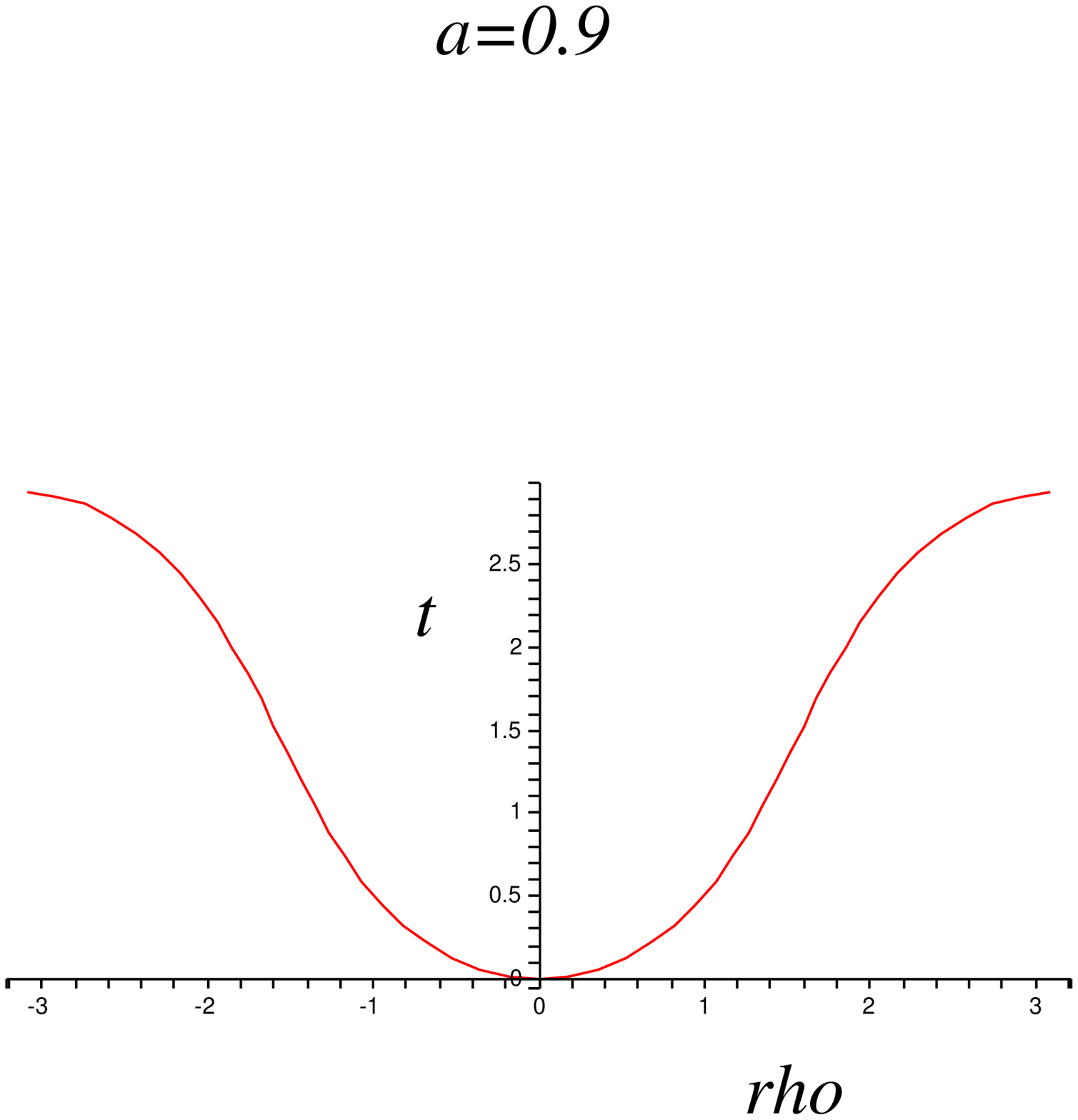}
\includegraphics[scale=0.45]{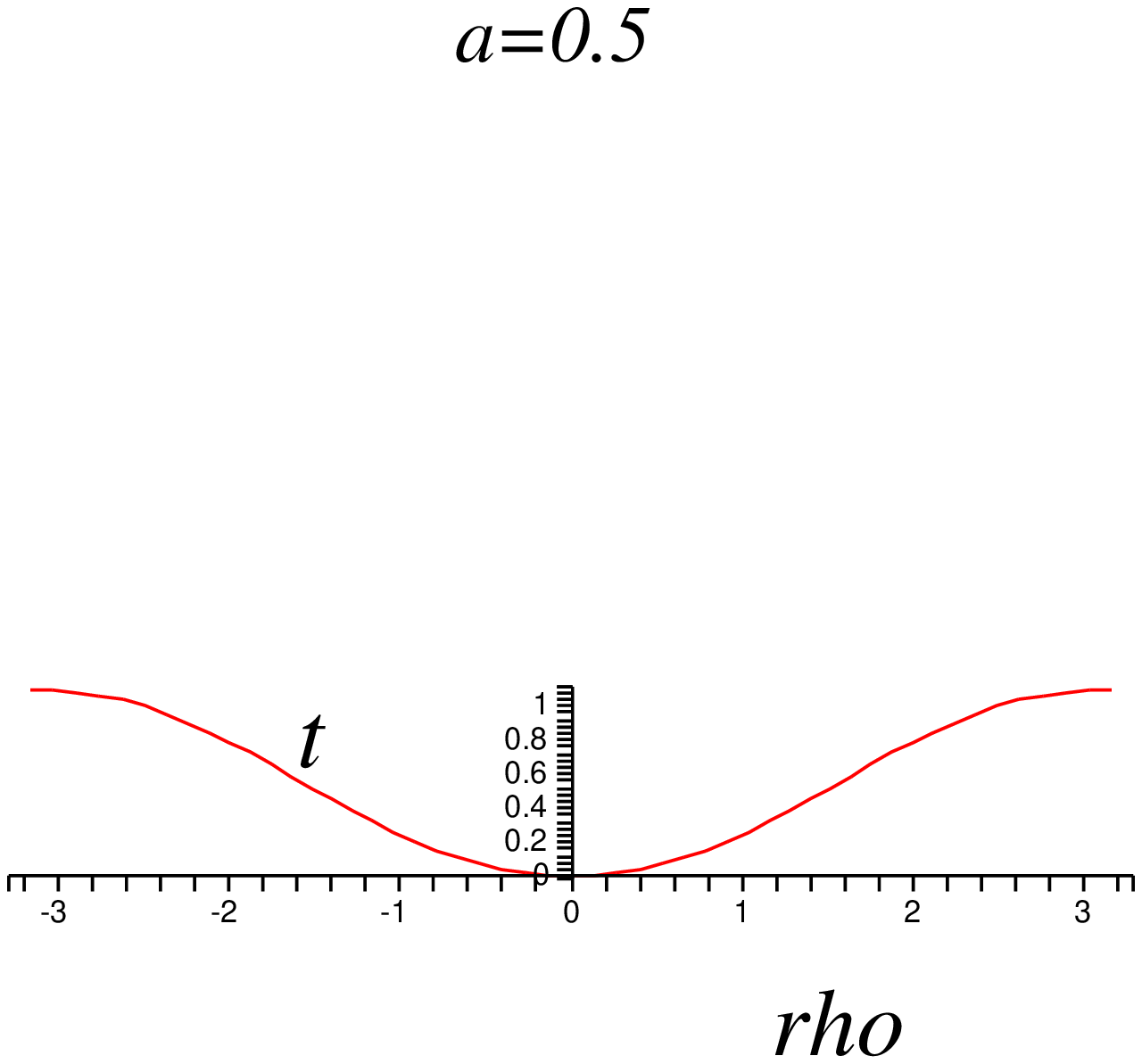}
\end{figure}

\Proof
With the notations of the lemma \ref{zero}, if $\delta <+\infty$ then,
since $\rho_a^\prime$ takes its values in $\mathopen]0,1\mathopen]$, 
$\rho_a$ would have a finite limit $l$ at $\delta$. Since 
$a\in \mathopen]0,1\mathopen[$, we have $1-a^2\sin^2 l>0$ which 
allows to extend $\rho_a$ beyond $\delta$, contradicting the
maximality of $\rho_a$. This shows that $\delta=+\infty$. 

Since $\rho_a^\prime \geq \sqrt{1-a^2}>0$ and $\rho_a (0)=0$, 
there is a smallest $s_1>0$ such that $\rho_a(s_1)=\pi$.

Now let us consider the function $f(s):=2\pi-\rho_a(2s_1-s)$, 
$s\in \R$. We observe that $f$ is also a solution of equation    
(\ref{intprem}) and satisfies
$f(s_1)=2\pi-\rho_a(s_1)=\pi=\rho_a(s_1)$. 
Consequently $f\equiv \rho_a$, that is:
\begin{equation}\label{rho1}
\rho_a(2s_1-s)=2\pi -\rho_a(s),\ \forall s \in \R.
\end{equation}
 As we are interested in curves generating rotational totally
umbilic surfaces, we look for a function $t$ satisfying 
$t^{\prime 2}= 1- \rho_a^{\prime 2}=a^2\sin^2  \rho_a$.
Let $t_a$ be the function defined on $\R$ by setting 
$t_a^\prime (s)=a\sin \rho_a(s)$ and $t_a(0)=0$.

Since $\rho_a$ is an odd function 
and $t_a(0)=0$ we deduce that $t_a$ is an even function.  
Observe that the function $g(s):=t_a(2s_1-s)$ satisfies
$g^\prime(s)=t_a^\prime(s)$ (using equation (\ref{rho1})) 
and $g(s_1)=t_a(s_1)$. Thus $g\equiv t_a$, that is
$t_a(2s_1-s)=t_a(s)$ for any $s\in \R$. Using the evenness of $t_a$
we get
\begin{equation}\label{t}
t_a(s+2s_1)=t_a(s),\ \forall s \in \R.
\end{equation}
Using equation (\ref{rho1}) and the oddness of $\rho_a$ we obtain
for any $s\in \R$
$$
\rho_a(s)=-\rho_a(-s)=-(2\pi -\rho_a(2s_1+s)),
$$
and so 

\begin{equation}\label{rho2}
\rho_a(s+2s_1)=2\pi +\rho_a(s),\ \forall s \in \R.
\end{equation}

Now the curve $\wt \alpha_a (s)=(\rho_a(s),t_a(s))$, $s \in \R$, is a curve in
the 
Riemannian universal cover $\wt \cm$ of $\cm$. Observe that the equations 
(\ref{t}) and (\ref{rho2}) show that restricting $s$ to 
$\mathopen[-s_1,s_1\mathopen]$, $\wt \alpha_a$ gives rise to an analytic closed
curve
$\alpha_a$ in $\cm$. Since $\rho_a$ is increasing on 
$\mathopen[-s_1,s_1\mathopen]$ and $\rho_a(-s_1)=-\rho_a(s_1)=-\pi$ we
deduce that $\alpha_a$ is embedded and nonhomologous to zero. 

As $\rho_a$ is odd and $t_a$ is even, the curve $\alpha_a$ has the
desired symmetry.

It is clear that the other choice 
$t_a^\prime (s)=-a\sin \rho_a(s)$ leads
to the curve deduced from $\alpha_a$ by the isometry 
$(\rho,t) \mapsto (\rho,-t)$. \qed

\begin{rem}\label{sym-geod}
We observe that the curve $\alpha_a$ is globally invariant under the isometry
$(\rho,t) \mapsto (\pi -\rho \ \mathrm{mod}(2\pi) , t_a(2s_0)-t)$.
\end{rem}

\begin{lem}\label{a=1}
Assume $a=1$. Then $\rho_1 (s)=\pi/2 - 2 \arctan e^{-s}$, $s \in \R$. 
This gives rise to a unique, up to an
ambient isometry, 
curve   $\alpha_1$ in $P\subset \cm$ generating 
a rotational totally umbilic surface.
The curve $\alpha_1$ is
complete, open, embedded and
symmetric with
respect to the axis of rotation $R$.
The rotational totally
umbilic surface, $S_1$, generated by $\alpha_1$ in 
$\sd \times \R$ is complete, properly embedded, analytic and 
homeomorphic to $\m 2$.
\end{lem}

\Proof
  Since $\rho_1$ is the maximal solution of 
$\rho^\prime  =\sqrt{1-\sin^2\rho}$ satisfying $\rho(0)=0$, we deduce that 
$\rho_1$ is solution of $\rho^\prime =\cos \rho$. A straightforward
  computation shows that the maximal solution of this last equation 
with the initial condition $\rho(0)=0$ is 
$$
\rho_1 (s)=\frac{\pi}{2} - 2 \arctan e^{-s},\ \forall s \in \R.
$$
Note that $\rho_1$ takes its values in 
$\mathopen]-\pi/2,\pi/2\mathopen[$.

As in the proof of the lemma \ref{a<1}, we can assume that 
$t$ satisfies $t^\prime (s) =\sin\rho_1(s)$, up to an ambient isometry. 
We consider the function $t_1$
defined by setting $t_1^\prime (s) =\sin\rho_1(s)$ and $t_1(0)=0$. It is
straightforward to check that $\sin\rho_1 (s)=\tanh (s)$ and then
$t_1(s)=\log \cosh (s)$. 

As $\rho_1$ is odd and $t_1$ is even, the curve $\alpha_1$ has the
desired symmetry. As a matter of fact, $\alpha_1$ is the graph of the
function $t(\rho)=-\log \cos \rho,\ \rho \in \mathopen]-\pi/2,\pi/2[$.
This concludes the proof. \qed

\begin{figure}[ht]
\includegraphics[scale=0.45]{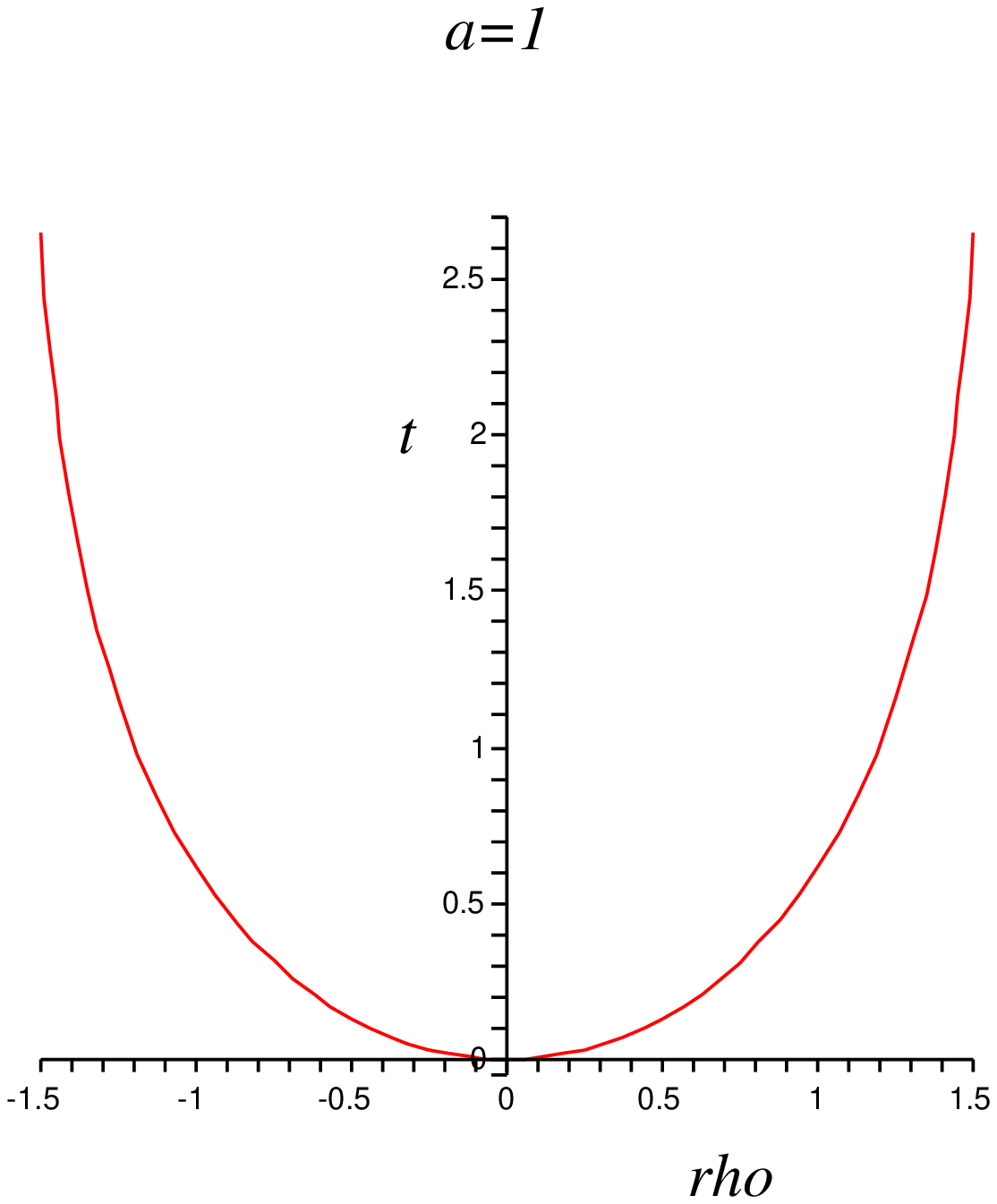}
\end{figure}

\begin{lem}\label{a>1}
Assume $a>1$. Then the maximal solution $\rho_a$ is defined on a
bounded interval $\mathopen]-\delta_a,\delta_a\mathopen[$ where
$\delta_a$ is a positive number. It gives rise to a unique, up to an
ambient isometry, 
Jordan curve   $\alpha_a$ in $P\subset \cm$ generating 
a rotational totally umbilic surface. The curve $\alpha_a$
is  analytic and symmetric with respect to
the axis of rotation $R$. 
The rotational totally
umbilic surface, $S_2(a)$, generated by $\alpha_a$ is 
analytic, embedded and homeomorphic to the sphere.
 Furthermore, $S_2(a)$ is homologous to
zero in $\sd \times \R$.
\end{lem}

\Proof
Since $a>1$ we deduce that 
$\rho_a(s) \in \mathopen]-\arcsin 1/a  ,\arcsin 1/a\mathopen[$.
Recall that $\rho_a$ is defined on an
open interval $\mathopen]-\delta_a,\delta_a\mathopen[$, 
see the lemma \ref{zero}. We first show that $\delta_a$ is finite.
Assume by contradiction $\delta_a=+\infty$. 
Since $\rho_a$ is
nondecreasing it admits a limit 
$l\in \mathopen] 0,\arcsin 1/a \mathopen]$
as $s\to +\infty$. Necessarily $l=\arcsin 1/a$ since otherwise 
it would follow from the equation (\ref{intprem}) that 
$\rho_a^\prime >\sqrt{1-a^2\sin^2l}  >0$ for all $s>0$. 
Therefore $\rho_a$ would not be bounded. 

Using the equation (\ref{second}) we see that for 
$s$ big enough : 
$\rho_a^{\prime \prime}(s)\leq -1/(2\tan l)<0$. 
Consequently $\rho_a^\prime$ would be negative for $s$ big enough
which is a contradiction. This proves that $\delta_a$ is finite.

Let us call again $l$ the limit of $\rho_a$ as $s \to \delta_a$. If 
$l< \arcsin 1/a$ then we could extend the solution $\rho_a$, which is
maximal, beyond $\delta_a$, which is absurd. So $l=\arcsin 1/a$ 
and $\rho_a^\prime \to 0$ as $s \to \delta_a$.

Observe that, since the function $\rho_a$ satisfies equations 
(\ref{second}) and (\ref{intprem}), it satisfies also the following
equation
\begin{equation}\label{simple}
\rho^{\prime \prime}=-a^2\sin \rho \cos \rho.
\end{equation}

As the second member of 
(\ref{simple}) is bounded, its maximal solutions are defined on the
whole of $\R$. Call $\wt \rho_a$ the maximal solution of
(\ref{simple}) extending $\rho_a$. Set 
$f(s):=\wt \rho_a (2\delta_a -s)$. It is clear that $f$ and $\wt \rho_a$
satisfy equation (\ref{simple}) with the same initial conditions at
$\delta_a$. Thus we have 
\begin{equation}\label{even}
\wt \rho_a (2\delta_a -s)=\wt \rho_a (s), \ \forall s \in \R.
\end{equation}

   As we are interested in curves generating rotational totally
   umbilic surfaces, we look for a function $t$ satisfying 
$t^{\prime 2}= 1-\wt \rho_a^{\prime 2}=a^2\sin^2 \wt \rho_a$.
Let $t_a$ be the function defined by $t_a^\prime =a\sin \wt \rho_a$
and $t_a(0)=0$. As $\wt \rho_a$ is an odd function we deduce that
$t_a$ is even. Observe that the function 
$g(s):=2t_a(\delta_a)-t_a(2\delta_a -s)$ satisfies 
$g^\prime (s)=t^\prime _a (s)$ (using equation (\ref{even})) and 
$g(\delta_a)=t_a(\delta_a)$. Thus $g\equiv t_a$, that is 
\begin{equation}\label{tood} 
t_a(2\delta_a-s)= 2t_a(\delta_a)-t_a(s),\ \forall s \in \R.
\end{equation}

It follows from (\ref{even}) and the oddness of $\wt \rho_a$ that 

\begin{equation*}
 \wt \rho_a (s+4\delta_a)=\wt \rho_a (s),\ \forall s \in \R.
\end{equation*}

In the same way, using (\ref{tood}) and the evenness of $t_a$, we get
\begin{equation*}
t_a(s+4\delta_a)=t_a(s),\ \forall s \in \R.
\end{equation*}

 Now the curve $\alpha_a(s)=(\wt \rho_a (s),t_a(s))$, $s\in \R$, 
parametrizes a closed analytic curve in $P\subset \cm$.
Taking into account (\ref{even}), the oddnes of $\wt \rho_a$, 
(\ref{tood}) and the evenness of $t_a$, we deduce that the curve
$\alpha_a$ is symmetric with respect to the axis $R$.

Considering the fact that $t_a$ is increasing on
$\mathopen[0,2\delta_a\mathopen]$ and the symmetry of $\alpha_a$,  
we infer that $\alpha_a$ defines a Jordan curve in $P$.

To conclude the proof observe that the other choice for $t_a$,
that is 
$t_a^\prime =-a\sin \wt \rho_a$, leads to the curve deduced 
from $\alpha_a$ by the isometry $(\rho ,t)\mapsto (\rho ,-t)$. \qed

\begin{figure}[ht]
\includegraphics[scale=0.50]{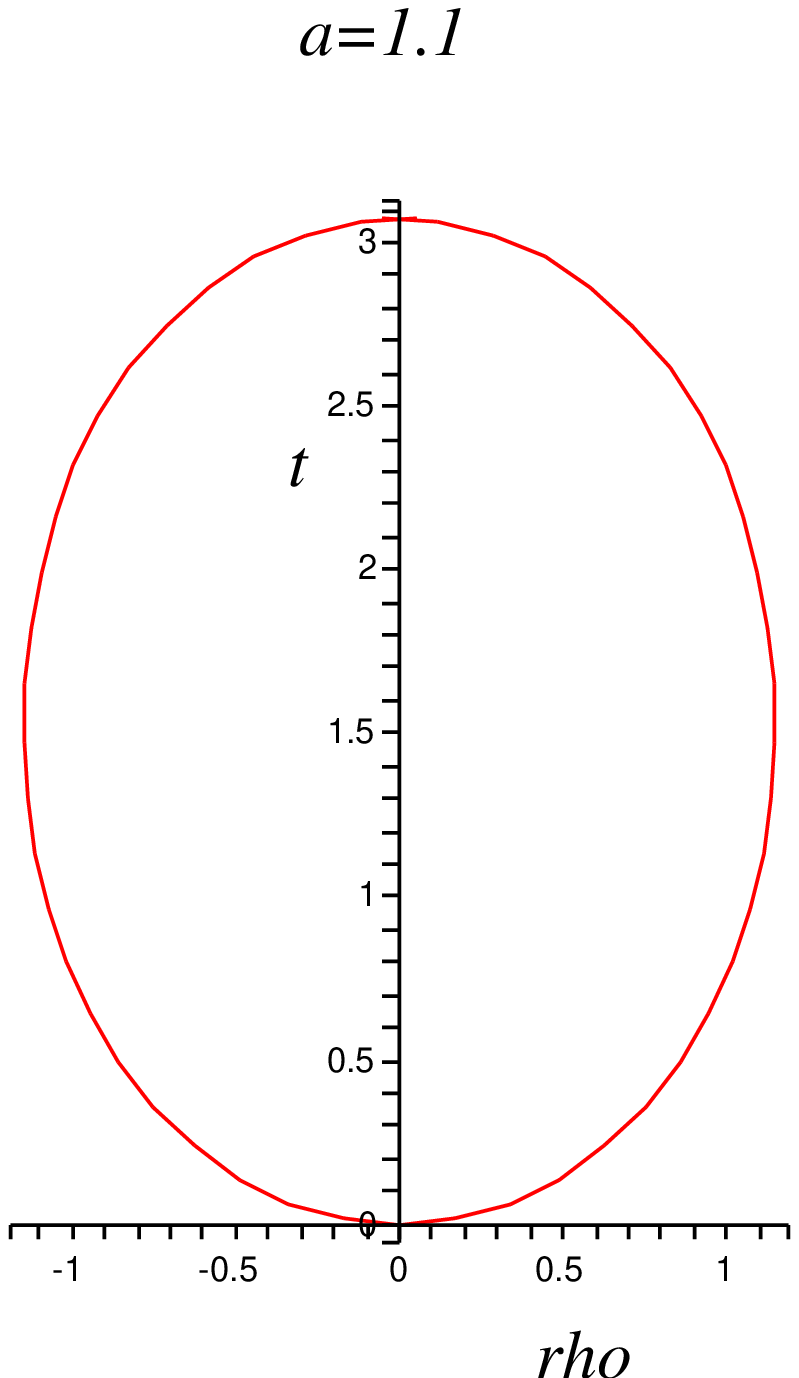} \qquad
\includegraphics[scale=0.50]{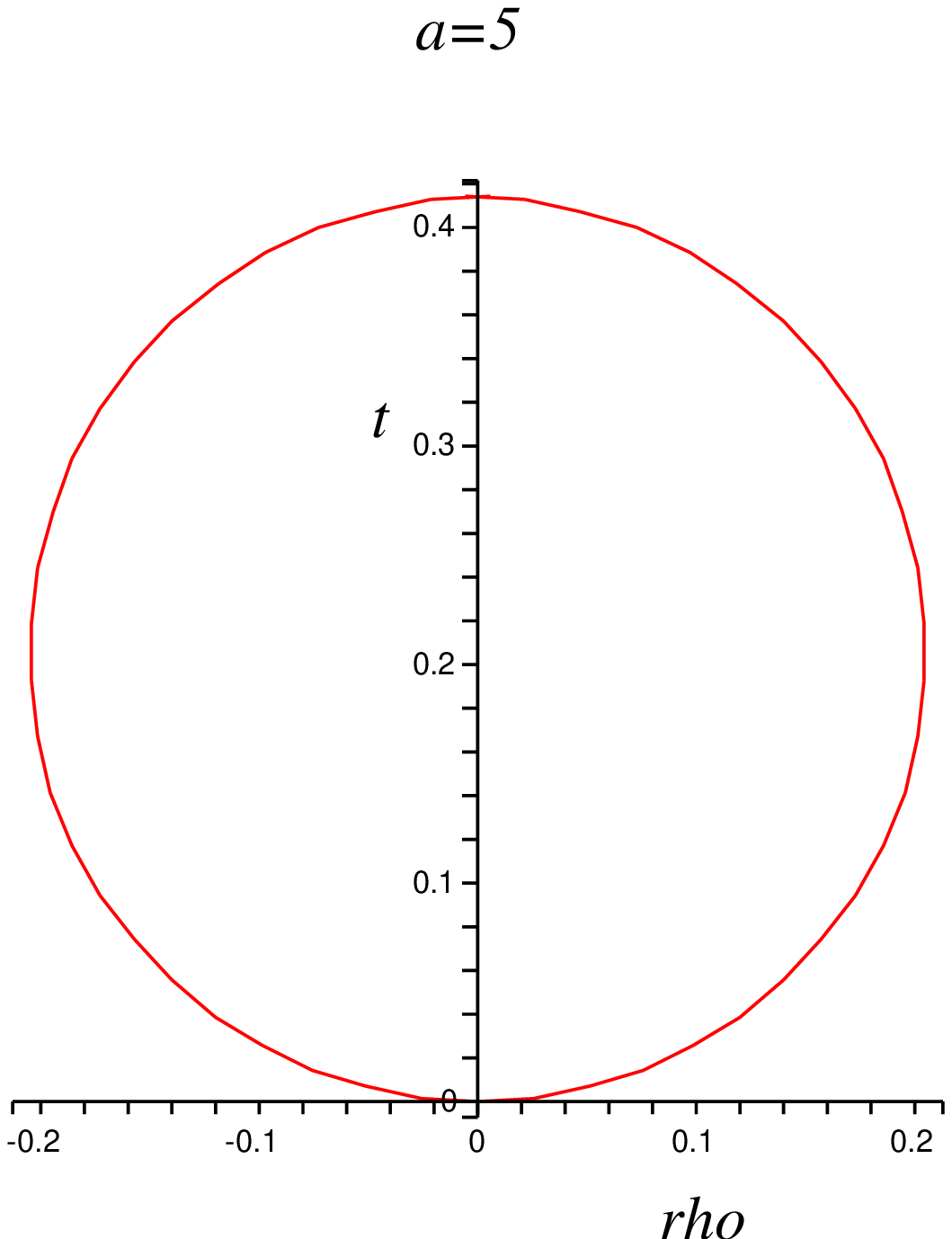}
\end{figure}

\begin{rem}\label{sym-plane}
The relations (\ref{even}) and (\ref{tood}) show that the curve
$\alpha_a$ in the lemma \ref{a>1} is symmetric with respect
to the horizontal reflection 
$(\rho,t)\mapsto (\rho, 2t_a(\delta_a) -t)$ in $P$.
Consequently, the surface $S_3(a)$ is symmetric with respect to the 
slice $\sd \times \{t_a (\delta_a)\}$.

 Furthermore, we observe that each surface $S_3(a)$ is contained 
in $\sd_{-}\times \R$ where $\sd_{-}$ is the south hemispshere.
\end{rem}

Summarizing we can state the following result.

\begin{thm}\label{S2}
Besides the totally geodesic slices $\sd \times \{t_0\}$ and the
vertical cylinder $\Gamma \times \R$ where $\Gamma \subset \sd$ is a
geodesic, 
the surfaces introduced in the lemmas \ref{a<1}, \ref{a=1} and 
\ref{a>1} are, up to ambient isometries, the only complete 
totally umbilic rotational surfaces in $\sd \times \R$. In particular
they are all embedded and homeomorphic either to $\m 2$ or to
$\sd$. Among the surfaces homeomorphic to $\sd$ some are
homologous to zero and some are not.
\end{thm}

\begin{rem}
It is interesting to observe that unlike in the case of space forms, 
the totally umbilic surfaces we obtained do not have constant mean
curvature, except for the totally geodesic ones.
\end{rem}

\section{Symmetric totally umbilic surfaces in $\hd (\kappa) \times \R$}
\label{plan hyperbolique}
In this section, we classify the totally umbilic surfaces in 
$\hd \times \R$ which are invariant under a one-parameter group of
isometries. The case of $\M^2(\kappa)\times \R$, $\kappa<0$, is
completely similar.

We recall that in $\hd$ there are three kinds of one-parameter
families
of positive isometries:
the rotations around a fixed point ({\sl elliptic isometries}),
the translations along a fixed geodesic ({\sl hyperbolic isometries})
and the "translations" along the horocycles sharing the same point at
infinity ({\sl parabolic isometries}).  An isometry of $\hd$ 
obviously induces an isometry of $\hd \times \R$ fixing the factor 
$\R$ pointwise. Such an isometry of $\hd\times \R$ obtained from 
an elliptic (resp. parabolic, hyperbolic) isometry of $\hd$
will thus be called elliptic (resp. parabolic, hyperbolic).

We will see that 
for each of the associated families of isometries of $\hd \times \R$
there exist complete and globally invariant
totally umbilic surfaces. In fact, we are going to
classify all of them. 
More precisely, we prove they are all 
embedded, those which are invariant under elliptic isometries are
either totally geodesic slices $\hd \times \{t_0\}$ or
homeomorphic to the sphere and the remaining ones are all homeomorphic
to the plane. In particular the only totally geodesic ones are the
slices and the products $\Gamma \times \R$ where $\Gamma \subset \hd$
is a geodesic.

 We will work with the disk model for $\hd $, so that 
$$
\hd =\{(x,y)\in \m 2,\ x^2+y^2<1 \},
$$
and the metric is 
$$
\rmd s\h ^2=\left( \dfrac{2}{1-(x^2+y^2)}\right)^2
(\rmd x^2 + \rmd y^2).
$$
Therefore the product metric on $\hd \times \R$ reads as follows:
$$
\rmd \tilde s ^2=\left( \dfrac{2}{1-(x^2+y^2)}\right)^2
(\rmd x^2 + \rmd y^2) + \rmd t^2,
$$
where $(x,y)\in \hd$ and $t\in \R$.
We consider the following particular geodesics of $\hd$ :
\begin{equation*}
\begin{aligned}
\Gamma = &\ \{(x,0), \ x\in \mathopen]-1,1\mathopen[\,\}\subset \hd \\
L=&\ \{(0,y), \ y\in \mathopen]-1,1\mathopen[\,\}\subset \hd \\
\end{aligned}
\end{equation*}

Up to ambient isometries, 
we can assume that the
symmetric surfaces 
are generated by curves in
the geodesic plane $P:= \Gamma \times \R \subset \hd\times \R$.

On the geodesic $\Gamma$ we denote by $\rho \in \R$ the signed
distance to the origin $(0,0)$, thus $x=\tanh \rho/2$. 
Therefore the metric on $P$ is 
$$
\rmd s^2=\left( \dfrac{2}{1-x^2}\right)^2\rmd x^2 +\rmd t^2
=\rmd\rho^2 +\rmd t^2.
$$

Given a curve
$\alpha(s)=(\rho (s),t(s))$ 
parametrized
by arclength in $P$, we let $\theta (s)$ be the oriented angle between
the $\rho$-axis and $\alpha^\prime (s)$.
Therefore, we have:
\begin{equation}\label{htangent}
\left\{ 
\begin{aligned}
\rho^\prime (s)=&\ \cos \theta(s)\\
t^\prime(s)=&\ \sin\theta(s)\\
\end{aligned}\right.
\end{equation}

In the elliptic case, the isometries of $\hd\times \R$
under consideration are the rotations around the vertical axis 
$R:= \{(0,0)\}\times \R$. In the parabolic case, the isometries are 
the ones corresponding to the point at infinity 
$(-1,0)\in \partial_\infty \hd$. Finally, the hyperbolic
isometries correspond to translations along $L$ in $\hd$.

In the plane $P$ we consider the unit normal $N$ to the curve 
$\alpha$ so that the basis ($\alpha^\prime (s),N(s))$ is positively oriented for
each $s$. In the three cases we orient by $N$ the symmetric
surface generated by $\alpha$. The principal curvatures 
computed with respect to this orientation are as follows:
$$
\lambda_1(s)=\ \theta^\prime (s) 
$$
and
$$
\lambda_2(s)=\ \left\{ 
\begin{matrix}
 \dfrac{\sin \theta (s)}{\tanh \rho (s)} &  \mathrm{(elliptic\ case)}\\
\sin \theta (s) &  \mathrm{(parabolic\ case)}\\
\sin \theta (s)\tanh \rho (s) & \mathrm{(hyperbolic\ case)}\\
\end{matrix}\right.
$$

\subsection{Elliptic case}{ $ $ }

The umbilicity condition  is 
\begin{equation*}
\theta^\prime (s) =\dfrac{\sin \theta (s)}{\tanh \rho (s)}\cdot
\end{equation*}

This case is similar to the case $a>1$ in $\sd \times \R$, so we will
omit the details.

Differentiating the first equation in (\ref{htangent}) and using the
umbilicity condition we get
\begin{equation}\label{Esecond}
\rho^{\prime \prime}=\dfrac{(\rho^{\prime 2}-1)}{\tanh \rho} \cdot
\end{equation}

Discarding the trivial totally geodesic surfaces
$\hd \times\{t_0\}$,
we can show as in the case of $\sd \times \R$
that $\rho^{\prime 2}(s)\not= 1$ for any $s$ such that 
$\rho (s) \not= 0$. Therefore we may assume that 
$\rho^{\prime 2}<1$. 

We can state the following.

\begin{prop}\label{P.Hyp.Elliptic}
Any local solution of (\ref{Esecond}) satisfying 
$\rho^{\prime 2}<1$ 
gives rise to a unique, up to 
ambient isometries, complete rotational totally umbilic 
and nongeodesic surface in 
$\hd \times \R$. Moreover, there exists a one-parameter family of such 
surfaces and all of them are analytic, embedded and homeomorphic to
the sphere.

Besides the totally geodesic slices $\hd\times \{t\}$, these surfaces
are the only complete rotational and totally umbilic surfaces in $\hd
\times \R$.

  Furthermore, any rotational umbilic (including geodesic) surface in
  $\hd \times \R$ is, 
up to an ambient isometry, part of one of the above surfaces.
\end{prop}

\Proof
Let $\rho$ be a local solution of  (\ref{Esecond}). 
Proceeding as in the case of $\sd \times \R$, we can suppose that
$\rho^\prime >0$ and so 

\begin{equation}\label{Eintprem}
\rho^\prime =\sqrt{1-b^2\sh^2\rho},
\end{equation}
for some real number $b>0$. As in the case of $\sd \times \R$,  the functions 
$\rho (s)$ and $t(s)$ are related to the Jacobi elliptic functions as
follows: 
$\rho (s)=-i\mathrm{am}(is,-b^2)$ and, up to the sign, 
$t^\prime(s)=ib\ \mathrm{sn}(is,-b^2)$, see 
\cite[Chapter 16]{[Abra-Ste]}. Again, we
prefer to give direct and elementary arguments.

Let $\rho_b$ be the maximal solution of 
(\ref{Eintprem}) extending $\rho$. As in the proof of the lemma
\ref{zero}, we can prove that $\rho_b$ vanishes at some point. Thus,
up to a reparametrization we can assume that $\rho_b(0)=0$,
consequently we prove as in lemma \ref{zero} that $\rho_b$ is an odd
function. Therefore $\rho_b$ is defined on an interval 
$\mathopen]-\delta_b,\delta_b\mathopen[$. As in lemma \ref{a>1}, it
can be shown that $\delta_b$ is a finite positive number, that
$\rho_b$ has a finite limit $l=\arg\!\sinh 1/b$ 
at $\delta_b$ and that
$\rho_b^\prime (s) \to 0$ as $s \to \delta_b$.  

\begin{figure}[ht]
\includegraphics[scale=0.45]{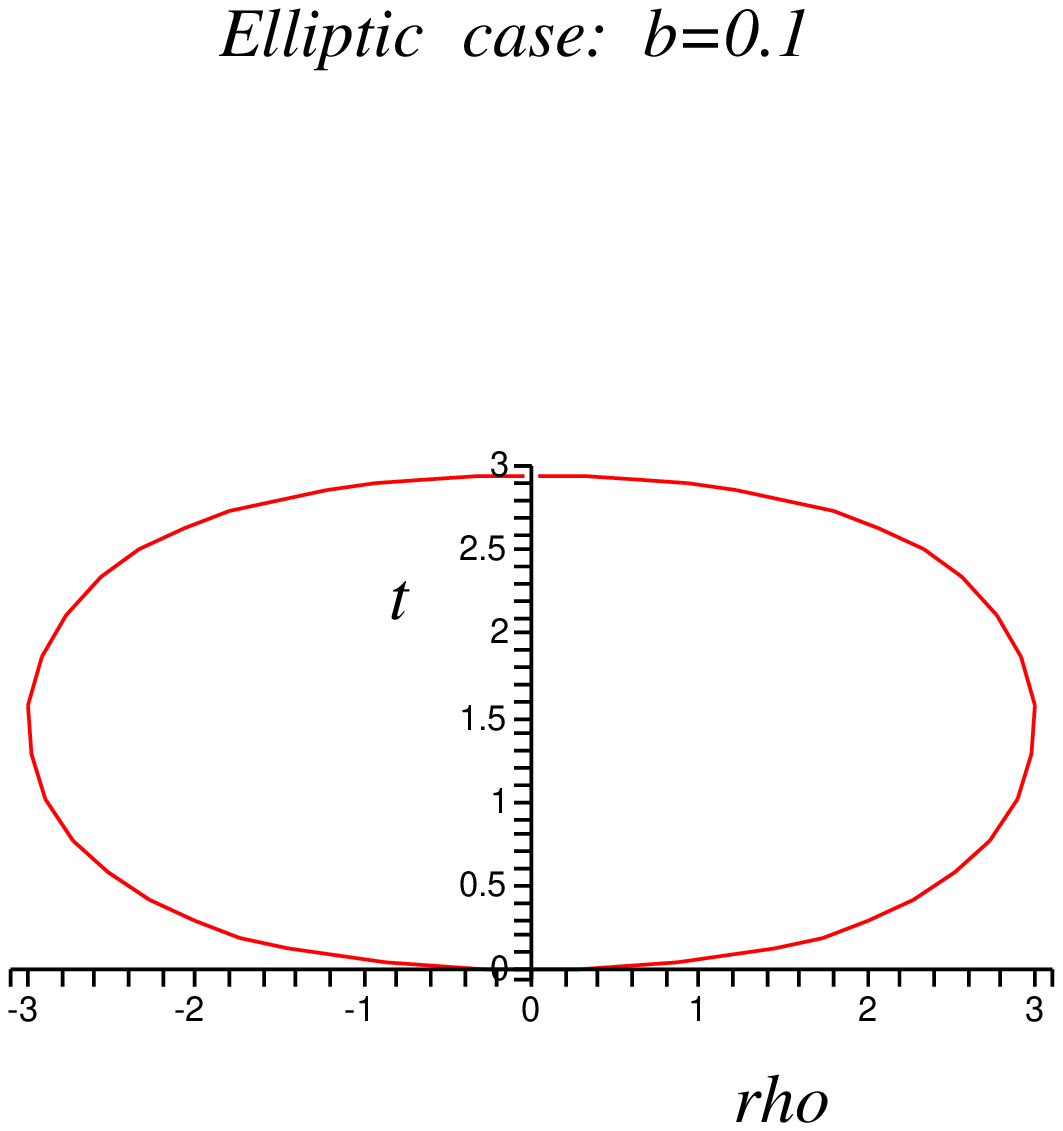} 
\includegraphics[scale=0.45]{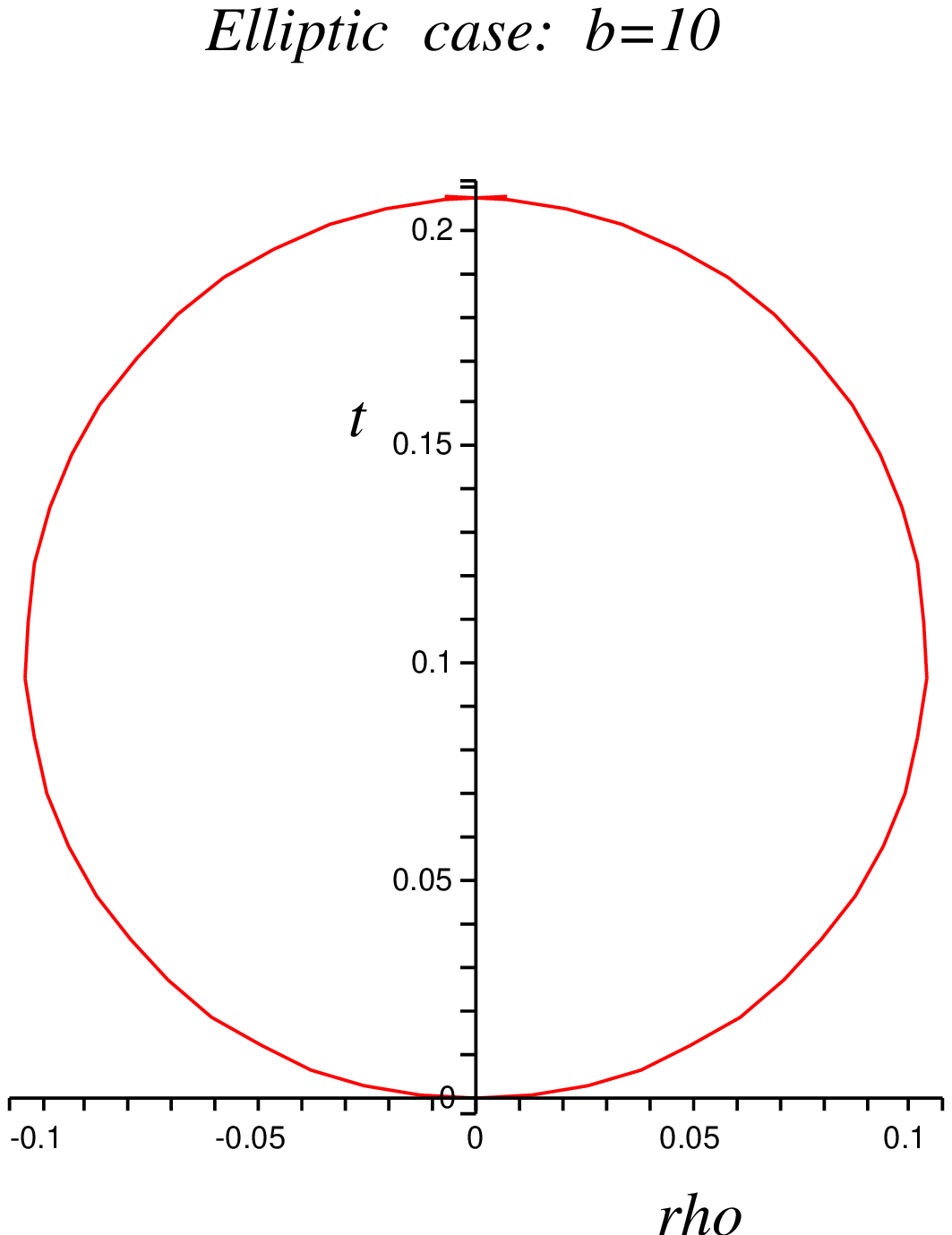}
\end{figure}


From equations (\ref{Esecond}) and (\ref{Eintprem}) we deduce that 
$\rho_b$ satisfies
\begin{equation}\label{Esimple}
\rho^{\prime \prime}=-b^2\ch \rho \sh \rho
\end{equation}
with the initial conditions $\rho (0)=0$ and $\rho^\prime (0)=1$.
Therefore we can extend the solution $\rho_b$ of (\ref{Esimple})
beyond $\delta_b$. Let $\wt \rho_b$ be the maximal solution of 
(\ref{Esimple}) extending $\rho_b$. Observe that for any $s_0$ where 
$\wt{\rho_b} ^\prime (s_0)=0$ we have the symmetry 
$\wt \rho_b (2s_0-s)=\wt \rho_b (s)$. As 
$\wt{ \rho_b}^\prime (\delta_b)=0$ and $\wt \rho_b$ is odd, 
we deduce that $\wt \rho_b$ is
defined on all of $\R$ and that it is $4\delta_b$-periodic.

As we are interested in curves generating rotational totally
   umbilic surfaces, we look for a function $t$ satisfying 
$t^{\prime 2}= 1-\wt{ \rho_b}^{\prime 2}=b^2\sh^2 \wt \rho_b$.
Let $t_b$ be the function defined by $t_b^\prime =b\sh \wt \rho_b$
and $t_b(0)=0$, thus $t_b$ is an even function.
As in the proof of lemma \ref{a>1}, we can show that $t_b$ 
satisfies $t_b(2\delta_b-s)=2t_b(\delta_b)-t_b(s)$ for any $s\in \R$ and that it is 
also $4\delta_b$-periodic.

Taking into account that $t_b$ is increasing on 
$\mathopen[ 0,2\delta_b\mathopen]$, we deduce that the curve 
$\alpha_b(s)=(\wt \rho_b (s),t(s))$, $s\in \R$, parametrizes an 
analytic Jordan curve in $P$, symmetric with respect to the axis $R$.

To conclude the proof we just observe that the other choice for $t_b$,
that is $t_b^\prime =-b\sh \wt \rho_b$, leads to the curve deduced 
from $\alpha_b$ by the isometry $(\rho ,t)\mapsto (\rho ,-t)$. \qed

\bigskip

\subsection{Parabolic case}{ $ $ }\label{parabol}

The umbilicity condition  is 
\begin{equation*}
\theta^\prime (s) =\sin \theta (s)
\end{equation*}

Integrating this equation we get
\begin{equation*}
\theta (s)= 2\arctan \lambda e^s, \ \forall s\in \R,
\end{equation*}
for some real constant $\lambda$. 

First observe that $\lambda=0$ leads to the curve $\Gamma$ 
which generates a slice $\hd \times \{t_0\}$. 

Now if $\lambda <0$ then $\theta$ is a negative function. Note that 
the symmetry $(\rho ,t)\mapsto (\rho, -t)$ changes $\theta$ into 
$-\theta$. Therefore, up to an ambient isometry, we can assume that 
$\theta$ is positive and then $\lambda >0$. Finally observe that up to
the reparametrization $s\mapsto s -\log \lambda$ we can assume
$\lambda=1$ and then:
\begin{equation}\label{Ptheta}
\theta (s)= 2\arctan  e^s, \ \forall s\in \R,
\end{equation}

Taking into account the first equation in (\ref{htangent}) and 
(\ref{Ptheta}) we obtain 
$$
\rho^\prime (s)=\cos 2\arctan e^s=-\tanh s,\ \forall s \in \R.
$$
Thus 
$$
\rho(s)=-\log \ch s +\mu,\ \forall s \in \R,
$$
for some real constant $\mu$. Note that the isometries
of $\hd \times \R$ obtained from the hyperbolic translations along 
the geodesic $\Gamma$ in $\hd$ send any surface invariant under
the parabolic isometries fixing the point at infinity 
$(-1,0)\in \partial_\infty \hd$ to a surface of the same type.
Consequently, up to an ambient isometry, we can assume $\mu=0$. Thus
$$
\rho(s)=-\log \ch s ,\ \forall s \in \R,
$$

As for the function $t$, taking into account the second equation 
in (\ref{htangent}) and 
(\ref{Ptheta}) we obtain 
$$
t^\prime (s)=\sin 2\arctan e^s=\dfrac{2e^s}{1+e^{2s}},
\ \forall s \in \R.
$$
Integrating we get 
$$
t(s)=2 \arctan e^s +\beta,\ \forall s \in \R,
$$
for some real constant $\beta$. Up to a vertical translation we can
take $\beta =-\pi/2$ so that $t(0)=0$ and 
$$
t(s)=2 \arctan e^s -\dfrac{\pi}{2},\ \forall s \in \R.
$$
Note that $t(-s)=-t(s)$ and $\rho (-s)=\rho (s)$ 
for any $s \in \R$ so that the curve $L_{par}$ parametrized by 
$(\rho (s),t(s)),\ s\in \R$, 
is symmetric with respect to $\Gamma$.



Summarizing we state the following proposition.

\begin{prop}
Besides the slices $\hd\times\{t\}$, up to ambient isometries,
there exists a unique complete 
totally umbilic surface $S_P$ in $\hd \times \R$ invariant under parabolic
isometries. This surface is analytic, properly embedded, homeomorphic to a
plane and is invariant under reflection with respect to a horizontal
slice.

Moreover, any totally umbilic surface invariant under parabolic
isometries is either part of a slice or, up to an ambient isometry, 
part of this surface.
\end{prop}

\begin{figure}[ht]
\includegraphics[scale=0.38]{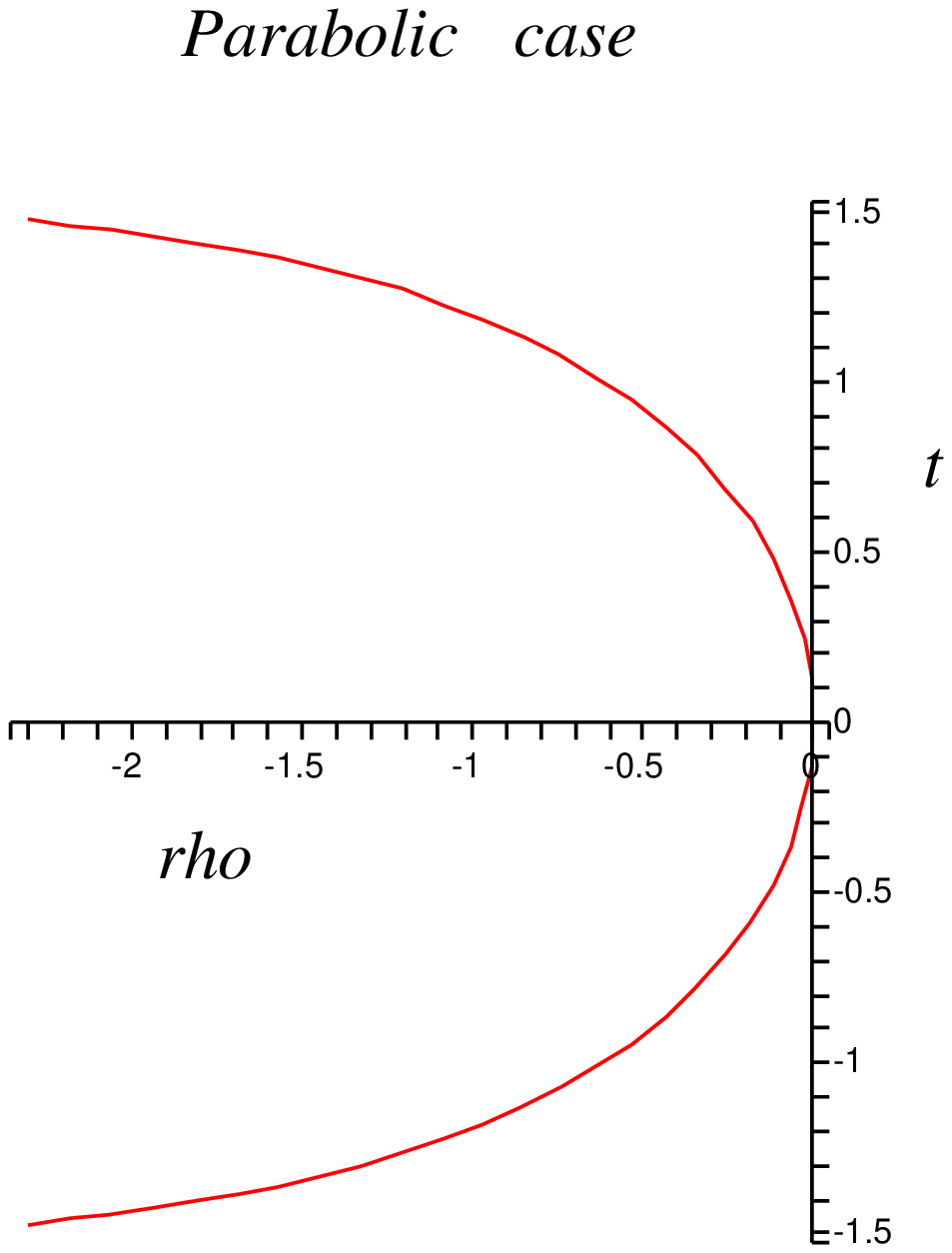}
\end{figure}

\begin{rem}
 Consider the surface $M$ in $\hd \times \R$ generated by the same
curve $L_{par}\subset P$ under parabolic isometries 
fixing, now, the point at infinity 
$(1,0)\in \partial_\infty \hd$ (and not $(-1,0)$ as before).
Observe that at each point of $M$ the
principal curvatures $\wt\lambda_1,\ \wt\lambda_2$ of $M$ are given by
$\wt\lambda_1=\lambda_1$ and $\wt\lambda_2=-\lambda_2$. Therefore we
get
$$
\wt\lambda_1 +\wt\lambda_2 =\lambda_1-\lambda_2 \equiv 0.
$$
We deduce that $M$ is a complete minimal surface of $\hd \times \R$, embedded and 
invariant under parabolic isometries. Consequently, $M$ is foliated by
horocycles.     
This minimal surface was
considered by B. Daniel \cite{[Daniel1]}, R. Sa Earp \cite{[Earp]} 
and L. Hauswirth \cite{[Hauswirth]}.
\end{rem}

\subsection{Hyperbolic case}{ $ $ }

The umbilicity condition is 
\begin{equation*}
\theta^\prime (s) =\sin \theta (s) \tanh \rho (s)
\end{equation*}

Proceeding as in the elliptic case, we discard the totally geodesic
surfaces $\hd \times \{t\}$ and therefore we can 
we assume that $\rho$ satisfies
$\rho^{\prime 2} <1$ and:
\begin{equation}\label{Eq.hyper}
\rho^{\prime \prime}=(\rho^{\prime 2}-1)\tanh \rho,\ \ 
\rho^{\prime 2}-1=-c^2\ch^2 \rho,
\end{equation}
for some real constant $c\in \mathopen]0,1[$. 
Again, invoking the Jacobi elliptic functions, it can be shown that, up to the
sign, we have  
$\rho (s)= i \mathrm{am}(is +K, c^2)-i\pi/2$ and 
$t^\prime (s)= c\ \mathrm{sn}(is +K, c^2)$ where 
$K=\int_0^{\pi/2} \frac{dt}{\sqrt{1-c^2\sin^2 t} }$, see 
\cite[Chapter 16]{[Abra-Ste]}.
Nevertheless, as in
the previuous cases, we prefer to give direct and elementary arguments.

We deduce from (\ref{Eq.hyper}) that 
\begin{equation}\label{Hsecond}
\rho^{\prime \prime}=-c^2\ch \rho \sh \rho.
\end{equation}
If $\rho \equiv 0$, then the generated surface is the vertical totally
geodesic plane \newline
$\{(0,y,t), \ -1<y<1,\  t\in \R\}$. Discarding this
case, we consider only the nontrivial solutions of (\ref{Hsecond}).

Again, as in the elliptic case, it can be shown that any maximal
solution of the last equation is defined on the whole of $\R$, is
periodic, vanishes at some point and, up to a reparametrization, is odd.
We can
therefore assume that there exists a unique maximal solution $\rho_c$
satisfying $\rho (0)=0$ and $\rho^\prime (0)=\sqrt{1-c^2}$.

As we are interested in curves generating totally
umbilic surfaces, we look for a function $t$ satisfying 
$t^{\prime 2}= 1- \rho_c^{\prime 2}=c^2\ch^2  \rho_c$.
Let $t_c$ be the function defined by $t_c^\prime =c\ch  \rho_c$
and $t_c(0)=0$. The function $t_c$ is odd. 
Moreover consider any $s_0 \in \R$ such that $\rho_c^\prime (s_0)=0$,
then it can be shown that $\rho_c(2s_0-s)=\rho_c(s)$.
Set
$T(s):=-t_c(2s_0-s)+2t_c(s_0)$, $s\in \R$. We have 
$T^\prime\equiv t_c^\prime$ and $T(s_0)=t_c(s_0)$, therefore
$T=t_c$. As the function $t_c$ is odd, we have 
$t_c(s+4s_0)=4t_c(s_0)+t_c(s)$ for every $s\in \R$.
We can deduce that $t_c$ is an increasing and
nonbounded function on $\R$. This shows that 
the curve $\alpha_c =(\rho_c,t_c)$ is properly embedded.

Observe that the other choice $t^\prime =- c\ch  \rho_c$ 
changes the curve $\alpha_c =(\rho_c,t_c)$ into 
the symmetric curve $(\rho_c,-t_c)$ with respect to $\Gamma$.

\begin{figure}[ht]
\includegraphics[scale=0.45]{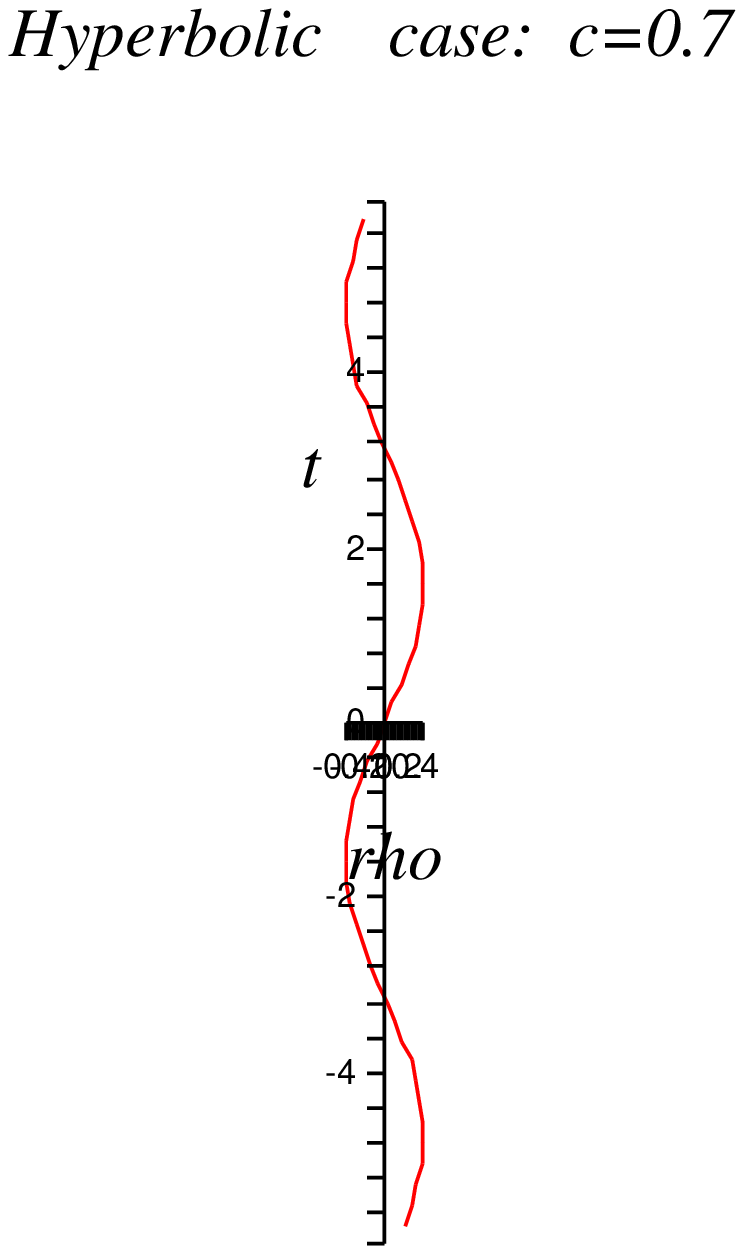}
\includegraphics[scale=0.45]{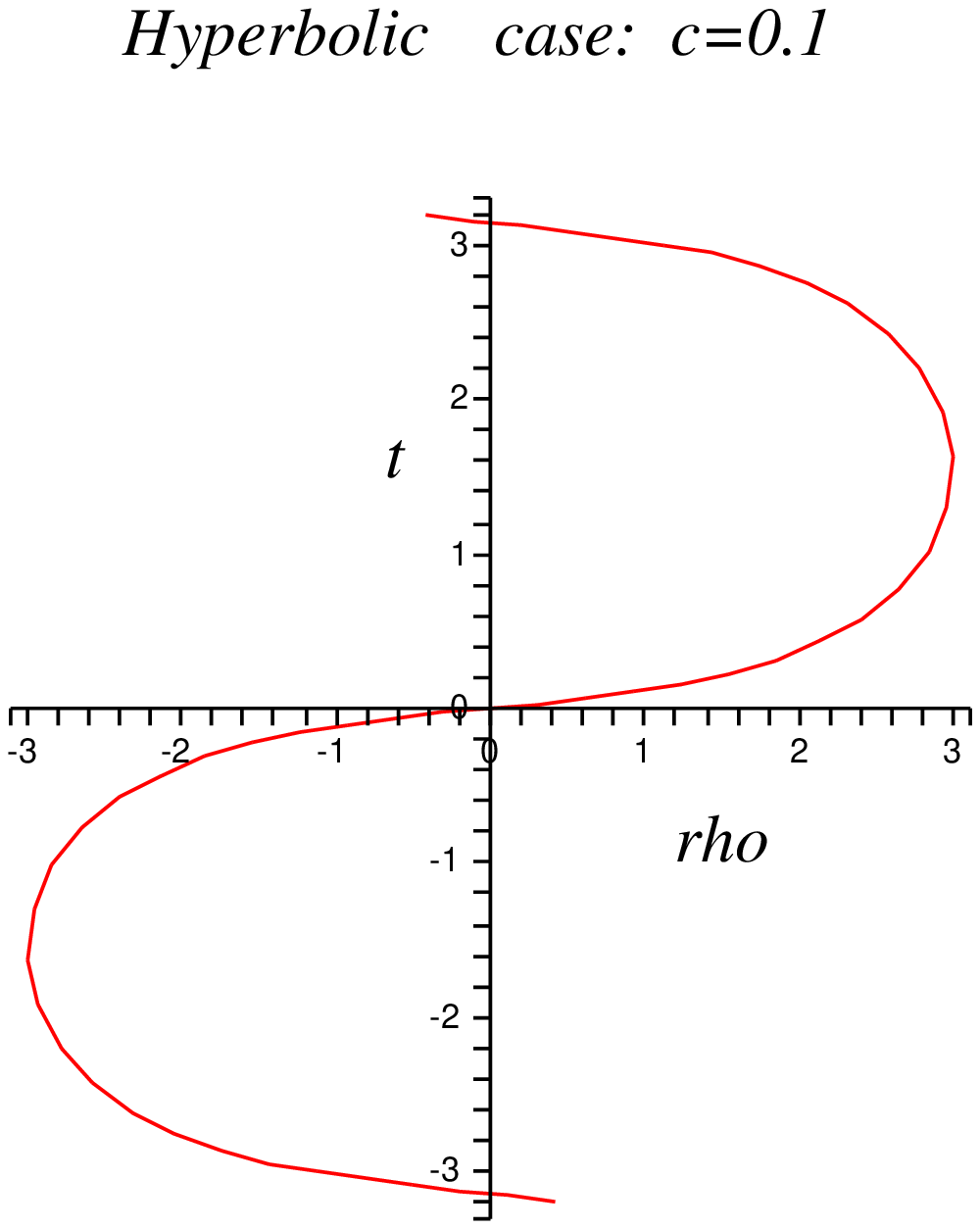}
\end{figure}


We summarize stating the following.

\begin{prop}
Any non identically zero 
local solution of (\ref{Hsecond}) satisfying 
$\rho^{\prime 2}<1$ 
gives rise to a unique, up to 
ambient isometries, complete totally umbilic and nongeodesic surface
$S_c,\ c>0$,  in 
$\hd \times \R$ invariant under hyperbolic isometries.
Moreover, there exists a one-parameter family of such 
surfaces and all of them are analytic, properly embedded 
and homeomorphic to
the plane. These surfaces are periodic in the vertical direction and
symmetric with respect to a discrete set of horizontal slices.

  Furthermore, any umbilic surface in $\hd \times \R$
invariant under hyperbolic isometries is either part of a vertical
totally geodesic plane, a slice $\hd \times \{t\}$ or, 
up to an ambient isometry, part of one of the surfaces $S_c$.
\end{prop}

\section{Unicity of totally umbilic surfaces in $\hd (\kappa) \times \R$ and  
$\sd (\kappa) \times \R$}\label{unicite}

In this section $\M^2$ stands for $\hd$ or $\sd$. The cases 
$\M^2=\M^2 (\kappa)$ for $\kappa <0$ and $\kappa >0$ are completely
analogous.

 We will need the
following result which is of independent interest.

\begin{prop}\label{independant}
Let $S\subset \M^2 \times \R$ be an orientable surface transversal to
each slice $\M^2\times \{t\}$. 
We suppose the following:
\begin{enumerate}
\item The geodesic curvature of each horizontal curve 
$S_t:= S \cap (\M^2\times \{t\})$ in $\M^2$ 
is constant (depending on $t$).
\item The angle between $S$ and $\M^2\times \{t\}$ is constant  
along $S_t$ for each $t$.
\end{enumerate}
Then:

In case $\M^2 =\sd$ the surface
 $S$ is part of rotational surface.

In case $\M^2 =\hd$ the surface $S$ is part of either 
\begin{enumerate}
\item a rotational surface,
\item or a surface invariant by a family of parabolic isometries 
having the same fixed point at infinity,
\item or a surface invariant by a family of hyperbolic isometries 
along the same fixed geodesic of $\hd$.
\end{enumerate}
\end{prop}

\Proof
Let $N$ be a unit normal field along $S$. We define the function $\nu$
on $S$ setting $\nu := \langle N,\pt \rangle$. We denote by $T$ the
projection of $\pt$ on $S$, that is $T=\pt-\nu N$. 

 As the angle between $S$ and $\M^2\times \{t\}$ is constant along
 $S_t$, we deduce that $S_t$ is a line of curvature. Indeed let 
$\gamma : s\in I\subset \R \rightarrow \gamma (s)\in S_t $ 
be a regular parametrization
 of $S_t$, then taking into account that $\pt$ is a parallel field and
 the definition of $T$:
\begin{equation*}
0=\frac{\rmd}{\rmd s}\langle N,\pt\rangle =
\langle \ov{\nabla}_{\gamma^\prime (s)}N,\pt \rangle=
\langle \ov{\nabla}_{\gamma^\prime (s)}N, T \rangle ,
\end{equation*} 
where $\ov{\nabla}$ is the connection on $\M^2 \times \R$.
It follows that $T$ is a principal direction on $S$.

Let $c: u\in I\subset \R \rightarrow c(u)\in S$ be a line of 
curvature associated
to the field $T$. We are going to show
that $c(I)$ is contained in a vertical totally geodesic plane. This is
equivalent to showing that the horizontal projection 
$c_h : I\subset \R \rightarrow \M^2 $ of $c$ is a geodesic.

Assume first that $c$ is never vertical, that is $\nu\not= 0$ along
$c$. Thus $c_h^\prime$ does not vanish. 

 Let $\nabla$ be the connection on $\M^2$. It is sufficient to show that 
$\nabla _{c_h^\prime}{c_h^\prime}$ is always parallel to
$c_h^\prime$. As $c^\prime=  c_h^\prime +(1-\nu^2)\pt$ we have
\begin{equation*}
\begin{aligned}
\ov{\nabla}_T T= &\ \ov{\nabla}_{c^\prime}c^\prime =\nabla _{c_h^\prime}{c_h^\prime} +
\frac{\rmd}{\rmd u}(1-\nu^2) \pt \\
=&\ \nabla _{c_h^\prime}{c_h^\prime} -2\nu \nu^\prime \pt\\
\end{aligned}
\end{equation*}
As $T$ is a principal direction there exists a function $\lambda$
such that $\ov{\nabla}_T N=\lambda T$. Therefore
$$
\nu^\prime=\frac{\rmd}{\rmd u} \langle N,\pt\rangle 
=\langle \ov{\nabla}_T N,\pt\rangle
=\lambda (1-\nu^2).
$$
Thus we obtain:
$$
\ov{\nabla}_T T=\nabla _{c_h^\prime}{c_h^\prime}
-2\lambda\nu (1-\nu^2)\pt \cdot
$$
Moreover we have 
\begin{equation*}
\begin{aligned}
\ov{\nabla}_T T= &\ \ov{\nabla}_T (\pt-\nu N) =-\ov{\nabla}_T \nu N\\
=&\ -\nu^\prime N-\lambda \nu T\\
=&\ (\frac{\nu\prime}{\nu}-\lambda \nu)c_h^\prime +
(-\nu^\prime \nu-\lambda \nu(1-\nu^2))\pt\\
=&\  (\frac{\nu\prime}{\nu}-\lambda \nu)c_h^\prime
-2\lambda \nu(1-\nu^2)\pt .
\end{aligned}
\end{equation*}
Thus we get
$$
\nabla _{c_h^\prime}{c_h^\prime}=
(\frac{\nu\prime}{\nu}-\lambda \nu)c_h^\prime ,
$$
which shows that $c_h(I)$ is a geodesic in $\M^2$.

We denote by $w$ a unit horizontal field along
$c$ tangent to $S$ and for each $u\in I$ 
we let $P(u)$ be the vertical totally geodesic
plane containing $c(u)$ and orthogonal at $c(u)$ to $w(u)$.

Suppose now that that $\nu$ vanishes on an open interval $J\subset I$.
Let $u_0\in J$. Observe that along the horizontal curve of $S$ through 
$c(u_0)$ the vector field $N$ is horizontal. This means that an open
set of $S$, including $c(J)$, is part of a cylinder $\gamma \times \R$
where $\gamma \subset \M^2$ is some horizontal curve. Clearly this
implies that $w$ is constant along $J$, and thus so is $P$.

Combining those two arguments we see that $P$ is locally
 constant on an open and dense subset of $I$. As $P(u)$ depends in a 
differentiable way on $u$, we conclude that $P$ is constant.

Let us now consider a horizontal curve $\gamma : I\rightarrow S_t$ 
parametrized by arclength. Let $s_1, s_2\in I$ and call 
$c: \mathopen]-\varepsilon,\varepsilon[ \rightarrow S$ 
the integral curve of $T$ such that $c(0)=\gamma(s_1)$ and 
$\wt c: \mathopen]-\varepsilon,\varepsilon[ \rightarrow S$ 
the integral curve of $T$ such that $\wt c(0)=\gamma(s_2)$.

Let us call $c_3$ (resp. $\wt c_3$) the vertical coordinate of $c$
(resp. $\wt c$). Calling again $u$ the parameter in 
$\mathopen]-\varepsilon,\varepsilon[$, we have
$$
c_3^\prime (u)=\langle c^\prime (u),\pt\rangle =\langle
T(c(u)),\pt\rangle
=1-\nu^2(c(u))=1-\nu^2(c_3(u)).
$$
Thus $c_3$ and $\wt c_3$ verify the same first order differential
equation with the same initial condition at $u=0$. We conclude that 
$c_3 \equiv \wt c_3 $.

Recall that $c$ and $\wt c$ are contained in vertical totally
  geodesic planes $P$ and $\wt P$.

Let us call $\Gamma \subset \M^2$ the complete constant geodesic
curvature line defined by $\gamma$, that is $\gamma \subset \Gamma$.

  Observe that there is a unique positive isometry $\varphi$ of $\M^2$ such that 
$\varphi (\Gamma)=\Gamma$, $\varphi (c(0))=\wt c(0)$ and preserving
the orientation of $\Gamma$.
 Therefore the
isometry \newline 
$\Phi (z,t)=(\varphi (z),t)$ of $\M^2 \times \R$ sends $P$ to
$\wt P$. Note that the curves $\wt c$ and $\Phi \circ c$ in the
vertical plane $\wt P$ have the same vertical component and make the
same angle with the horizontal for each 
$u \in \mathopen]-\varepsilon,\varepsilon[$. We deduce that these 
curves coincide: $\Phi \circ c=\wt c$. This concludes the proof. \qed

\medskip
We now state the main result of this section.

\begin{thm}\label{caracterisation}
Let $S\subset \M^2 \times \R$ be an immersed totally umbilic surface.
Then $S$ is part of a complete and embedded totally umbilic surface
$\wt S$ 
which is invariant by a one-parameter group of isometries of 
$\M^2 \times \R$. More precisely, up to an ambient isometry, in case
 $\M^2=\sd$, then 
$\wt S$ is one of the examples described
in the section \ref{sphere}, and in case $\M^2=\hd$ then 
$\wt S$ is one of the examples described
in the section \ref{plan hyperbolique}.

 In particular, any totally geodesic surface is part of 
a slice $\M^2\times \R$ or part of a product $\Gamma \times \R$ where
$\Gamma\subset \M^2$ is a geodesic.
\end{thm}

\Proof 
Locally $S$ is the image of an embedding 
$X:\Omega \rightarrow \M^2 \times \R$, where $\Omega$ is an open disk
in $\m 2$.
 As $X$ is
totally umbilic, there exists a function 
$\lambda :\Omega \rightarrow \R$ such that 
\begin{equation*}
\ov{\nabla}_w N=\lambda w,
\end{equation*}
for any vector $w$ tangent to $S$.

Proceeding as in the proof of Theorem \ref{nonexistence}, 
as $\tau =0$ we obtain:
\begin{equation}\label{gradient1}
\nabla \lambda =\kappa \nu T,
\end{equation}
where $\kappa$ is the Gaussian curvature of $\M^2$, that is $\kappa =1$
if $\M^2=\sd$ and $\kappa =-1$ if $\M^2=\hd$.

Assume for the moment that $\lambda$ has no critical point. 
In particular each level curve of $\lambda$ is orthogonal to $T$ and
is therefore horizontal, that is 
belongs to some $\M^2 \times \{t_0\}$.
Let $\gamma :I \subset \M^2 \times \{t_0\}$ be such a curve
parametrized by arclength. We have
\begin{equation*}
\frac{\rmd \nu}{\rmd s}=
\langle \ov{\nabla}_{\gamma^\prime (s)}N,\pt\rangle +
\langle N,\ov{\nabla}_{\gamma^\prime (s)}\pt \rangle 
=\lambda (\gamma(s)) \langle \gamma^\prime (s), \pt \rangle=0,
\end{equation*}
therefore $\nu$ is constant along $\gamma$.

We now call $n$ the unit normal field  
in $T(\M^2 \times \{t_0\})$ along $\gamma$ with the orientation induced
by $N$. Let $\theta$ be the oriented angle between $n$ and $N$,
hence 
$N(\gamma(s))=\cos \theta (s) n(\gamma(s)) +
\sin \theta (s)\pt(\gamma(s))$,  
we deduce that $\theta$ is
constant along $\gamma$.

On the other hand:

\begin{equation*}
\begin{aligned}
\lambda (\gamma(s)) =&\ 
\langle \ov{\nabla}_{\gamma^\prime  (s)}N,\gamma^\prime (s)\rangle\\
=&\ \langle \ov{\nabla}_{\gamma^\prime  (s)}
 (\cos \theta  \  n +\sin \theta \pt), \gamma^\prime  (s) \rangle\\
=&\ \cos \theta 
\langle \ov{\nabla}_{\gamma^\prime  (s)} n,\gamma^\prime  (s)\rangle 
+ \sin \theta \langle \ov{\nabla}_{\gamma^\prime  (s)}\pt,
\gamma^\prime  (s)\rangle\\
=&\  \cos \theta 
\langle \ov{\nabla}_{\gamma^\prime  (s)} n,\gamma^\prime  (s)\rangle.
\end{aligned}
\end{equation*}

Now observe that 
$\langle \ov{\nabla}_{\gamma^\prime  (s)}n,\gamma^\prime  (s)\rangle$ 
is the geodesic curvature of $\gamma$ in  $\M^2 \times \{t_0\}$. Since
$\lambda$ and $\theta$ are constant along $\gamma$ we deduce that
$\gamma$ has constant geodesic curvature. We conclude using the proposition 
\ref{independant} and results in the sections \ref{sphere} and
\ref{plan hyperbolique} that $S$ is as stated.

 Suppose now that $\lambda$ has some critical points. 

Let $U\subset S$ be a connected component, if any, of the interior of 
the set of critical points of $\lambda$. The formula (\ref{gradient1}) shows that
 $N$ is either always vertical or always horizontal in $U$. In the former
 case $U$ is part of a slice $\M^2 \times \{t_0\}$ and in the latter
 case $U$ is part of a cylinder, that is part of a product 
$\Gamma \times \R$ where $\Gamma$ is some curve in $\M^2$. As $S$ is
 totally umbilic, $\Gamma$ has to be a geodesic and so $U$ is totally
 geodesic.

Let now $V\subset S$ be a connected component, if any, of the set 
of regular points  of $\lambda$.  From the
first part of the proof, we know that $V$ is part of one the 
symmetric examples given in the sections \ref{sphere} and 
\ref{plan hyperbolique}.

Therefore $S$ is obtained by gluing pieces of totally geodesic
surfaces and pieces of the symmetric examples constructed in the
sections \ref{sphere} and \ref{plan hyperbolique}.
A closer look at these
different types of surfaces shows that the whole of $S$ is either 
totally geodesic or part of one of the complete symmetric examples, 
which concludes the proof. \qed

\begin{rem}
The local existence of totally umbilic surfaces in $\sd \times\R$ and 
$\hd \times \R$ can be seen in an alternate way. 
Indeed, it is known that umbilicity is
preserved by conformal diffeomorphisms, see \cite{[Spivak]}
(Vol. 4). It can be
  shown that $\sd \times \R$ is conformally diffeomorphic to 
$\R^3 \setminus \{(0,0,0)\}$, see Section \ref{conform}. This implies the umbilic
surfaces in 
$\sd \times \R$ correspond through this conformal diffeomorphism to those of 
$\R^3 \setminus \{(0,0,0)\}$. However to classify them in $\sd \times \R$ up to 
congruences, and to understand their geometry requires a nontrivial work. 
Regarding $\hd \times \R$, it can be shown that 
$\hd \times \mathopen]0,\pi[$ is 
conformally diffeomorphic to $\hi3$, see Section \ref{conform}. Nonetheless this is
not enough
to understand the global geometry and topology of the umbilic surfaces in $\hd
\times \R$.
\end{rem}

\section{Totally umbilic surfaces in $Sol$}\label{Sol}
The $Sol$ geometry is the eighth model geometry of Thurston, see 
\cite{[Thurston]}. It is a Lie group endowed with a left-invariant
metric, it is a homogeneous simply connected \newline
3-manifold with a
3-dimensional isometry group, see \cite{[Bonahon]}. It is 
isometric to $\R^3$ equipped with the metric:
$$
ds^2=e^{2z}\rmd x^2 +e^{-2z}\rmd y^2 +\rmd z^2.
$$
The group structure of $Sol$ is given by
$$
(x^\prime,y^\prime,z^\prime)\star (x,y,z)=
(e^{-z^\prime}x + x^\prime, e^{z^\prime}y + y^\prime, z+z^\prime).
$$
The isometries are:
$$
(x,y,z)\mapsto (\pm e^{-c}x +a ,\pm e^{c}y + b,  z+c)\ \ 
\textrm{and}\ \ 
(x,y,z)\mapsto (\pm e^{-c}y +a ,\pm e^{c}x + b,  -z+c),
$$
where $a,\ b$ and $c$ are any real numbers.
We set $E_1=e^{-z}\px$, $E_2=e^{z}\py$ and $E_3=\pz$.
Thus $(E_1,E_2,E_3)$ is a global orthonormal frame. A straightforward
computation gives:

\begin{equation}\label{derivees}
\begin{matrix}
\ov{\nabla}_{E_1}E_1 =\ -E_3 &  &
\ov{\nabla}_{E_2}E_1 =\ 0 &  &  \ov{\nabla}_{E_3}E_1 =\ 0\\
\ov{\nabla}_{E_1}E_2 =\ 0 &  &
\ov{\nabla}_{E_2}E_2 =\ E_3 & &  \ov{\nabla}_{E_3}E_2 =\ 0\\
\ov{\nabla}_{E_1}E_3 =\ E_1 &  &
\ov{\nabla}_{E_2}E_3 =\ -E_2 & &  \ov{\nabla}_{E_3}E_3 =\ 0\\
\end{matrix}
\end{equation}

We deduce from (\ref{derivees}) that the vertical planes $\{x=x_0\}$
and $\{y=y_0\}$ are totally geodesic complete surfaces and that the
horizontal planes $\{z=z_0\}$ are not totally umbilic surfaces (in
fact they are minimal surfaces). 

We now look for totally umbilic surfaces which are invariant under the
one parameter group of isometries 
$(x,y,z) \mapsto (x+c,y,z)$. Clearly, such a surface is generated by a
curve $\gamma$ in the totally geodesic plane $\{x=0\}$. 
Discarding the trivial case of a vertical plane $\{y=y_0\}$,  we can
assume that $\gamma$ locally is a graph over the $y$-axis. Thus $\gamma$ is
given by $\gamma (y)=(0,y,z(y))$. Therefore the generated surface is 
parametrized by:
\begin{equation*}
X(t,y):= (t,y,z(y)).
\end{equation*}
We have $X_t=(1,0,0)=e^z E_1$ and 
$X_y=(0,1,z^\prime)=e^{-z} E_2 +z^\prime E_3$. 
As a unit normal field we can take 
$$
N= \frac{e^z z^\prime}{\sqrt{1+e^{2z} {z^\prime}^2}}E_2
-\frac{1}{\sqrt{1+e^{2z} {z^\prime}^2}}E_3.
$$
We have:
\begin{equation*}
\begin{aligned}
\ov{\nabla}_{X_t} N=&\ -\frac{1}{\sqrt{1+e^{2z} {z^\prime}^2}} X_t\\
\ov{\nabla}_{X_y} N=&\ 
\frac{e^{-z}}{(1+e^{2z} {z^\prime}^2)^{3/2}}
(1+2e^{2z}{z^\prime}^2 +e^{2z} z^{\prime \prime})E_2 \\
&\ \hskip8mm
+\frac{z^\prime}{(1+e^{2z} {z^\prime}^2)^{3/2}}
(1+2e^{2z}{z^\prime}^2 +e^{2z} z^{\prime \prime})E_3
\end{aligned}
\end{equation*}
So that $X$ is a totally umbilic immersion if and only if 
$$
\ov{\nabla}_{X_y} N=  -\frac{1}{\sqrt{1+e^{2z} {z^\prime}^2}} X_y,
$$
that is if and only if 
\begin{equation}\label{Solumbilic1} 
z^{\prime \prime} +3 {z^\prime}^2+ 2e^{-2z}=0.
\end{equation}
A first integral of (\ref{Solumbilic1}) is 
\begin{equation*}
{z^\prime}^2=a e^{-6z}-e^{-2z},
\end{equation*}
where $a$ is any positive real number.

Assume $z^\prime (y_0)=0$ for some $y_0$. Considering the function 
$f(y)=z(2y_0 -y)$, we can see that the curve
$\gamma$ is symmetric with respect to the vertical geodesic
$\{y=y_0\}$. Therefore, up to the isometry $(x,y,z)\mapsto (x,-y,z)$ and 
restricting the domain of $z$ if needed, we can assume $z^\prime >0$. 
Therefore $z$ satisfies 
\begin{equation}\label{Solumbilic2} 
z^\prime=e^{-z}\sqrt{ae^{-4z}-1}.
\end{equation}
We consider the maximal solution of (\ref{Solumbilic2}) defined by $z$
and we call it again $z$, it is defined on an open interval 
$\mathopen]y_1,y_2[,\ -\infty \leq y_1 < y_2 \leq +\infty$. By
(\ref{Solumbilic2}) the function $z$ is 
bounded above and is increasing we deduce using (\ref{Solumbilic1})
that 
$z^{\prime  \prime}$ has a negative limit at $y_2$. Taking into
account the fact
that $z^\prime$ is positive, we deduce that $y_2$ is finite, 
$y_2<+\infty$. Moreover, since $z^\prime$ is a positive and decreasing
function, it has a nonnegative limit at $y_2$. If this limit were
positive, we could extend the solution $z$ of (\ref{Solumbilic2})
beyond $y_2$ which contradicts the maximality of $z$. Thus we have
$\lim_{y\to y_2}z^\prime (y)=0$ and consequently 
 $\lim_{y\to y_2}z (y)=\frac{1}{4} \log a$.

Consider now the
maximal solution  of (\ref{Solumbilic1})
defined by $z$ and call it $z_a$. As 
$z_a^\prime (y_2)=0$ we have $z_a(2y_2-y)=z_a(y)$. Up to the horizontal
translation $(x,y,z)\mapsto (x,y-y_2,z)$, which is an ambient isometry,
we can assume that $y_2=0$ and therefore $z_a$ is an even fonction and
is defined on an interval $\mathopen]-y_a,y_a[$ where
$0<y_a\leq +\infty$. Observe that there exist $A>0$ and $y_3>0$ such
that $z_a^{\prime \prime}(y)<-A<0$ for any $y>y_3$. Therefore, if 
$y_a=+\infty$ we have $\lim_{y\to +\infty}z_a(y)=-\infty$. Suppose now
that $y_a$ is finite. If $z_a$ had a finite limit at $y_a$ then
$z_a^\prime$ would have also a finite limite but then we could extend
the solution $z_a$ beyond $y_a$, which is absurd. We deduce that in both
cases, that is $y_a<+\infty$ and $y_a=+\infty$, we have 
$\lim_{y\to y_a}z_a(y)=-\infty$.

  Now we show that $y_a<+\infty$. Indeed, as $z_a$ is a solution of 
(\ref{Solumbilic1}) satisfying $z_a(0)=0$ we have 
$$
z^\prime (y)=-\sqrt{a}e^{-3z}\sqrt{1-\frac{e^{4z}}{a}},
$$
for any $y>0$. Since $\lim_{y\to y_a}z_a(y)=-\infty$, we deduce
that for some $y_4>0$ we have
$$
\frac{1}{\sqrt{a}}z^\prime e^{3z} < -\frac{1}{2}
$$
for any $y>y_4$. Therefore we have 
$$
\frac{1}{3\sqrt{a}}e^{3z}<-\frac{y}{2}+c
$$
for some real constant $c$ and for any $y>y_4$. This implies that
$y_a<+\infty$.

 Call $\gamma_a$ the graph of the function 
$z_a:\gamma_a:=\{(0,y,z_a(y)),\ -y_a<y<~y_a\}$ and call $F_a$ the totally
umbilic complete surface generated by $\gamma_a$:
$$
F_a:=\{(x,y,z_a(y)),\ x\in \R,\ -y_a<y<y_a\}
$$

  Let $z_a$ and $z_b$ be two maximal solutions of (\ref{Solumbilic1})
where $a$ and $b$ are any real numbers. Set 
$c=\frac{1}{4}\log\frac{b}{a}$ and consider the ambient isometry
$(x,y,z)\mapsto (e^{-c}x,e^{c}y,z+c)$. Observe that this isometry 
maps the planar curve $\gamma_a$ onto the planar curve 
$\gamma_b$ and maps any Euclidean line parallel to the $x$-axis onto a
line of the same type. Therefore this isometry maps the totally umbilic
surface $F_a$ onto the totally
umbilic surface $F_b$.

Summarizing, we state the following result.

\begin{prop}\label{Sol1}
Up to ambient isometries, there exist only two complete totally
umbilic surfaces in the $Sol$ group invariant under the horizontal
translations $(x,y,z)\mapsto (x+t,y,z),\ t\in \R$. The first one is
the totally geodesic plane $\{y=0\}$. The second one is nongeodesic,
is contained in a slab delimited by two totally geodesic planes 
$\{y=\pm y_0\},\ y_0>0$, and is asymptotic to these planes. Moreover it
is symmetric with respect to the totally geodesic plane $\{y=0\}$.
\end{prop}

\begin{figure}[ht]
\includegraphics[scale=0.45]{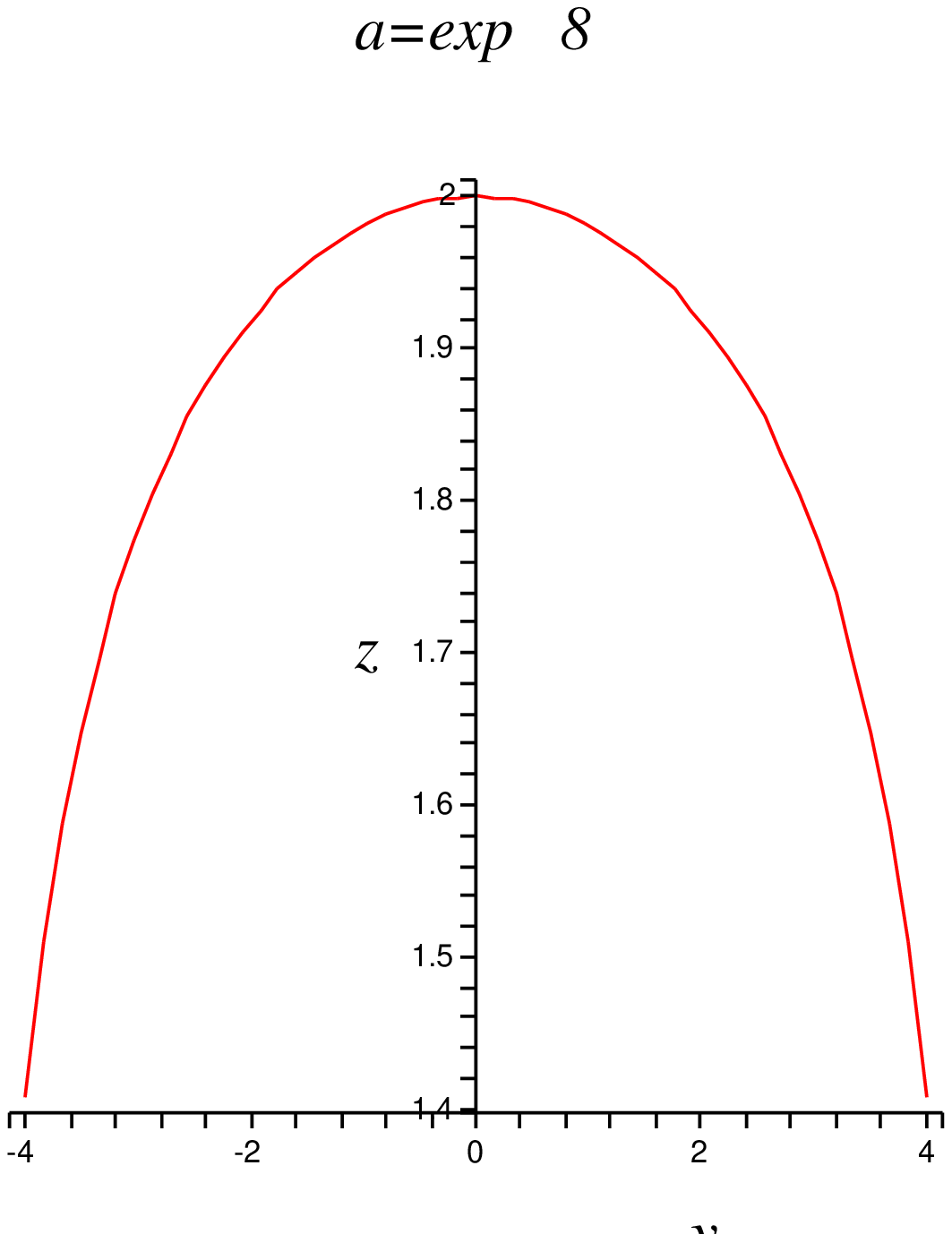}
\end{figure}


As a matter of fact we have the following.

\newpage

\begin{thm}
Up to ambient isometries, any totally umbilic surface in the $Sol$ group
is part of one of the two complete totally umbilic surfaces given in
the proposition \ref{Sol1}. In particular, up to ambient isometries,
there exists a unique complete totally geodesic surface in the $Sol$ group.
\end{thm}

\Proof
 Let $S$ be an immersed totally umbilic surface in 
the $Sol$ group.
Locally $S$ is the image of an embedding 
$X:\Omega \rightarrow Sol$, where $\Omega$ is an open disk
in $\m 2$. Call $(u,v)$ the coordinates on $\Omega$ and consider
a unit normal field $N$ on $X(\Omega)$. As $X$ is
totally umbilic, there exists a function 
$\lambda :\Omega \rightarrow \R$ such that 
\begin{equation*}
\left\{ 
\begin{aligned}
\ov{\nabla}_{X_u}N=\lambda X_u \\
\ov{\nabla}_{X_v}N=\lambda X_v
\end{aligned} \right.
\end{equation*}
where $\ov{\nabla}$ is the Riemannian connection of the $Sol$ group. As
in the proof of the theorem \ref{nonexistence} we find 

\begin{equation}\label{Soltenseur}
R(X_u,X_v)N=\lambda_v X_u - \lambda_u X_v,
\end{equation}
where $R$ denotes the curvature tensor of the $Sol$ group.
Let us express the later. let $X,Y,Z$ and $W$ be any vector fields.
Proceeding as in \cite{[Daniel]}, Proposition 2.1, after some
computations and using (\ref{derivees})
we obtain the
following:
\begin{equation*}
\begin{aligned}
\langle R(X,Y)Z,W\rangle =&\ 
(\langle X, Z \rangle \langle Y, W \rangle -
\langle X, W \rangle \langle Y, Z \rangle )\\
&\ +2\bigg(\langle X, W \rangle \langle Y, \pz \rangle \langle Z, \pz \rangle
+\langle Y, Z \rangle \langle X, \pz \rangle \langle W, \pz \rangle\\
&\ - \langle X, Z \rangle \langle Y, \pz \rangle \langle W, \pz \rangle
- \langle Y, W \rangle \langle X, \pz \rangle \langle Z, \pz \rangle \bigg).
\end{aligned}
\end{equation*}
We define the function $\nu$
on $\Omega$ setting $\nu := \langle N,\pz \rangle$. We denote 
by $T$ the projection of $\pz$ on $S$, that is $T=\pz-\nu N$.
We then have:
\begin{equation*}
\begin{aligned}
R(X_u,X_v)N=&\ 2\nu(\langle X_v,\pz\rangle X_u-
\langle X_u,\pz\rangle X_v) \\
=&\ 2\nu(\langle X_v,T\rangle X_u-
\langle X_u,T\rangle X_v)
\end{aligned}
\end{equation*}
From what we deduce using (\ref{Soltenseur}):
\begin{equation}\label{Solgradient}
\nabla \lambda =2\nu T.
\end{equation}

Assume first that $\nu$ and $T$ do not vanish on $\Omega$. Thus $T$ is of the
form 
$$
T=\alpha E_1 +\beta E_2 +\gamma E_3,
$$
where $\alpha$ and $\beta$ do not vanish simultaneously. 
Since $\vert T\vert ^2=1-\nu^2$ we have  
$\alpha^2+\beta^2=\nu^2(1-\nu^2)$. 
We thus have 
$$
  N= -\frac{\alpha}{\nu}E_1 -\frac{\beta}{\nu}E_2 +\nu E_3 .
$$
We set 
\begin{equation*}
 JT= -\frac{\beta}{\nu}E_1 +\frac{\alpha}{\nu}E_2
\end{equation*}
therefore $JT$ is tangent to $S$, orthogonal to $T$ and $E_3$ and 
satisfies $\vert JT\vert ^2=\vert T\vert ^2$. Furthermore we have 
\begin{equation*}
N\wedge T=JT,\ T\wedge JT=(1-\nu^2)N,\ JT\wedge N=T.
\end{equation*}
We now compute the derivative $[T,JT](\lambda)$ in two 
different ways.

We first compute $[T,JT]=\ov{\nabla}_T JT-\ov{\nabla}_{JT} T$. 
We have
\begin{equation*}
\begin{aligned}
\ov{\nabla}_T JT=&\ \ov{\nabla}_T N \wedge T +
N\wedge \ov{\nabla}_T T\\
=&\ N\wedge \ov{\nabla}_T T\ \ 
(\mathrm{since}\  \ov{\nabla}_T N =\lambda N)\\
=&\ N\wedge\ov{\nabla}_T (E_3 -\nu N)\\
=&\ N\wedge \ov{\nabla}_T E_3-\lambda \nu JT\\
\end{aligned}
\end{equation*}
Furthermore, using (\ref{derivees}) we obtain
\begin{equation*}
\begin{aligned}
\ov{\nabla}_T E_3 =&\ \alpha \ov{\nabla}_{E_1} E_3
+ \beta \ov{\nabla}_{E_2} E_3
+(1-\nu^2)\ov{\nabla}_{E_3} E_3 \\
=&\ \alpha E_1-\beta E_2\\
\end{aligned}
\end{equation*}
from what we deduce after some straightforward computations
$$
\ov{\nabla}_T E_3=\frac{\alpha^2-\beta^2}{1-\nu^2}T
-2\frac{\alpha \beta}{\nu(1-\nu^2)}JT -\nu(1-\nu^2)N.
$$
Consequently:
$$
\ov{\nabla}_T JT=2\frac{\alpha \beta}{\nu(1-\nu^2)} T
+\frac{\alpha^2-\beta^2}{1-\nu^2} JT-\lambda \nu JT.
$$
In the same way we obtain:
\begin{equation*}
\begin{aligned}
\ov{\nabla}_{JT} T=&\ \ov{\nabla}_{JT} E_3- \ov{\nabla}_{JT} \nu N \\
=&\ -2\frac{\alpha \beta}{\nu(1-\nu^2)}T +\nu JT
-\lambda \nu JT\\
\end{aligned}
\end{equation*}
We deduce that 
\begin{equation*}
[T,JT]=4\frac{\alpha \beta}{\nu(1-\nu^2)}T
+(\frac{\alpha^2-\beta^2}{1-\nu^2}-\nu)JT.
\end{equation*}
Using this last expression and (\ref{Solgradient}), we find
\begin{equation}\label{Lielambda1}
[T,JT] (\lambda) =8\alpha \beta .
\end{equation}

On the other hand, using again (\ref{Solgradient}) we have

\begin{equation}\label{Lielambda2}
\begin{aligned}
  \  [T,JT] (\lambda)=&\ T( JT(\lambda)) -JT( T(\lambda))\\
=&\ -JT( \langle \nabla \lambda, T\rangle)\\
=&\ 2(-1+3\nu^2)JT(\nu)\\
=&\ 4\frac{\alpha\beta}{\nu}(-1+3\nu^2)\\
\end{aligned}
\end{equation}
since an easy computation gives $JT(\nu)=2\alpha\beta/\nu$.
From (\ref{Lielambda1}) and (\ref{Lielambda2}) we deduce
$$
\alpha\beta(3\nu^2-2\nu-1)=0.
$$
Observe that if $\nu$ is constant on an open subset then 
$JT(\nu)\equiv 0$ which implies $\alpha \beta \equiv 0$.
So in all cases we have $\alpha \beta \equiv 0$.

Recall that $\alpha$ and $\beta$ do not vanish 
simultaneously since by our assumption $\nu\not=0$.
Therefore we have either $\alpha\equiv 0$ or 
$\beta\equiv 0$.

Considering the isometry 
$(x,y,z)\mapsto (y,x,-z)$ we can 
assume that $\alpha\equiv 0$. Then the surface is
part of a product $\R\times \Gamma$ where $\Gamma$ is a curve in the
geodesic plane $\{x=0\}$. This case is considered in the proposition
(\ref{Sol1}). 

Let us suppose now that $T\equiv 0$ on an open set. Then this open set
is part of a horizontal plane $\{z=z_0\}$, but this contradicts the
assumption of umbilicity.

To finish the proof we consider the case where $\nu\equiv 0$ on an
open subset. Therefore $T\equiv E_3$ and so this piece of the surface 
is  part of a product $L\times \R$ where $L$ is a curve in the 
horizontal plane $\{z=0\}$. If $L$ is contained in a line parallel 
to the $y$-axis, then the surface is contained in a totally geodesic
plane $\{x=x_0\}$. Discarding this trivial case, we can assume that 
$L$ is a graph over the $x$-axis. Consequently, the embedding $X$ is given 
by 
$$
X(x,t)=(x,y(x),t).
$$
As a unit normal we take 
\begin{equation*}
\begin{aligned}
N=&\ \frac{1}{\sqrt{e^{-2t}y^{\prime 2}+e^{2t}}}
(e^{-2t}y^\prime,-e^{2t},0)\\
=&\ \frac{y^\prime}{\sqrt{y^{\prime 2}+e^{4t}}}E_1
-\frac{1}{\sqrt{e^{-4t}y^{\prime 2}+1}}E_2
\end{aligned}
\end{equation*}
As $X_t=E_3$, using (\ref{derivees}) we obtain:
$$
\ov{\nabla}_{X_t}N=
-2\frac{y^\prime e^{4t} }{(y^{\prime 2}+e^{4t})^{3/2}}E_1
-2\frac{e^{-4t}y^{\prime 2}}{(e^{-4t}y^{\prime 2}+1)^{3/2}}E_2.
$$
The condition $\ov{\nabla}_{X_t}N =\lambda X_t$ is therefore
equivalent to $\lambda \equiv 0$ and $y^\prime \equiv 0$. So
$L$ is part of a line parallel to the $x$-axis and the surface is part
of a geodesic plane $\{y=y_0\}$. This concludes the proof. \qed

\begin{rem}\label{regularite}
It can be
proved that any twice differentiable
totally umbilic surface in a space form,
in $\sd (\kappa) \times \R$ or in $\hd(\kappa) \times \R$,
is in fact $C^3$ and then analytic by the previous discussions, see
\cite{[S-T]}.
\end{rem}

\section{An application}\label{conform}

As an application of the classification of totally umbilic surfaces obtained in
the previous sections, 
we can prove  the following result:

\begin{thm}\label{conformal}
Any conformal  diffeomorphism of $\hd\times\R,$ $\sd\times\R$ and Sol is an isometry.
\end{thm}

\Proof The result for $\sd\times\R$ is a consequence of the fact that the mapping:
\begin{equation*}
\begin{aligned}
\sd \times \R &
\rightarrow 
\m 3 -\{(0,0,0)\} \\
(p,t) & \mapsto  e^t p
\end{aligned}
\end{equation*}
is a conformal diffeomorphism -here $\sd$ is viewed as the unit sphere of  $\R^
3$ 
centered at the origin.
 Indeed the conformal diffeomorphisms of $\R^3 -\{(0,0,0)\}$ are the M\"obius
transformations fixing 
$(0,0,0)$ or sending $(0,0,0)$ to the point at infinity and these
transformations  correspond through the above conformal diffeomorphism to
isometries of $\sd\times\R.$ We leave the details to the reader.

We now prove the result for the space  $\hd\times\R.$  We claim that, except for
the slices 
$\hd\times\{t_0\},$ all the  non-compact maximal (for the inclusion) totally
umbilic surfaces  in $\hd\times\R$ are conformal to $\C.$ This is clear for the
products 
$\gamma\times\R,$ where $\gamma\subset \hd$ is a geodesic. As for the surfaces
invariant under a 
one parameter group of parabolic transformations and which are all congruent to
the surface $S_P$ described in \ref{parabol}, this is seen as follows. Consider
in $\mathbb{H}^3$ a totally geodesic 
plane which we call $\hd$ and denote by $N$ a unit normal along it. We let exp
denote the 
exponential map in $\mathbb{H}^3.$ 
Then the map: 
\begin{equation*}
\begin{aligned}
\hd \times \mathopen]0,\pi \mathopen[ &
\rightarrow 
\mathbb{H}^3\\
(p,t) & \mapsto  \text{exp}_p (\ln(\tan (\frac{t}{2})) N(p))
\end{aligned}
\end{equation*}
is a conformal diffeomorphism (cf. \cite{[S-T]} for the details) which sends
$S_P$ onto 
a totally umbilic surface of $\mathbb{H}^3$ with one point at infinity, that is
a horosphere. It remains to consider the case of the surfaces invariant under a
one parameter group of hyperbolic transformations. Consider such a surface
$\Sigma.$ We know that $\Sigma$ is 
invariant under a set of reflections of $\hd\times\R$ through horizontal slices 
$\hd\times\{t=a+nb\},$ for all $n\in \Z,$ and 
$a, b\in \R$ depending  on $\Sigma.$ Suppose by contradiction that $\Sigma$ 
is conformal to $\hd.$ The isometries of $\Sigma$
induced by those reflections  correspond then to conformal diffeomorphisms of
$\hd$ 
and so to isometries of $\hd.$ In particular all the horizontal curves
$\Sigma\cap \{ t=an+b\}$ correspond to geodesics of $\hd.$ Now observe that all
these curves are invariant by the hyperbolic isometries leaving $\Sigma$
invariant. We thus get isometries of $\hd$ 
which leave (globally) invariant more than one geodesic. This is a
contradiction 
as only the identity of $\hd$ has this property.

Take now a conformal diffeomorphism $\phi$ of $\hd\times\R.$ Then $\phi$ sends
any 
 horizontal slice $\hd\times\{t_0\} $
to a maximal totally umbilic non-compact surface which is conformal to $\hd.$
 From what preceeds it follows that $\phi$ sends  $\hd\times\{t_0\} $
conformally to 
some horizontal slice  $\hd\times\{t_1\} $
 and so isometrically. As $\phi$ is conformal this implies that for any 
$x\in \hd\times\{t_0\},$ the tangent map $D_x\phi$ is an isometry. So $\phi$ is
an isometry of $\hd\times\R.$ 
 
 The case of the {\it Sol}  group is treated analogously. There are two 
maximal 
totally umbilic surfaces  up to congruences. With the notations of Section
\ref{Sol},  the first one is the totally geodesic plane $\{ x=x_0\}$ and is
easily seen to be isometric to 
 $\mathbb {H}^2.$ The second one is the surface parametrized by:  $ X(t,y)= (t,y,z(y)),$
 where $t\in \R,\quad y\in \mathopen]-y_1,y_1\mathopen[$ and $z$ is the maximal 
solution to the equation: $z^\prime =e^{-z}\sqrt{e^{-4z}-1}$ (we have chosen
$a=1$ with the notations of Section \ref{Sol}).  The metric on this surface
writes: 
 $$ds^2 = e^{2z} dt^2+ e^{-z} dy^2.$$
 Making the change of coordinate $\xi =\int e^{-4z(y)}dy,$ the metric writes:
  $$ds^2 =e^{2z} (dt^2+d\xi^2).$$
 As  the function $z$ is bounded from above (cf. Section \ref{Sol}), the flat 
metric $dt^2+d\xi^2$ is complete. It follows that 
 the surface is conformal to $\C.$ As the totally geodesic planes $\{x=
const.\}$  
fill the whole space {\it Sol}, we conclude as before that any conformal
diffeomorphism of {\it Sol} is an isometry.  \qed


\begin{thebibliography}{10000}


\bibitem{[Abra-Ste]} M. Abramowitz and  A. Stegun.
{\em Handbook of mathematical functions}, Dover Publications, Inc., New York,
1972.

\bibitem{[Abresch-R]} U. Abresch and H. Rosenberg.
{\em A Hopf differential for constant mean curvature surfaces in 
$\sd \times \R$ and $\hd \times \R$}, Acta Math. 193, N$^o$ 2,
141-174, 2004.


\bibitem{[Bonahon]} F. Bonahon. {\em Geometric structures on
    3-manifolds}, Handbook of geometric topology, 93-164,
North-Holland, Amsterdam, 2002.


\bibitem{[Galvez]} J. Aledo, J. Espinar and J. G{\'a}lvez. 
{\em Complete surfaces of constant curvature in 
$\hd \times \R$ and $\sd \times \R$},  Calc. Var. Partial Differential Equations
 29,  no. 3, 347--363, 2007.


\bibitem{[Rato-CMC]} R. Caddeo, P. Piu and A. Ratto.
{\em SO(2)-invariant minimal and constant mean curvature surfaces
in 3-dimensional homogeneous spaces}, Manuscripta Math, 87, 1-12,
1995. 




\bibitem{[Daniel1]} B. Daniel. {\em Isometric immersions into 
$\Sp^n \times \R$ and $\hi n \times \R$ and applications 
to minimal surfaces}, Preprint, 2004.


\bibitem{[Daniel]} B. Daniel. {\em Isometric immersions into 
3-dimensional homogeneous manifolds},  Commentarii Math. Helv.  82, 
no. 1, 87--131, 2007.




\bibitem{[Earp]} R. Sa Earp, {\em Parabolic and hyperbolic screw
    motion surfaces in $\hi2 \times \R$}, 
http://www.mat.puc-rio.br/~earp/pscrew.pdf, 
to appear in J. Australian 
Math. Soc.


\bibitem{[Earp-T]}  R. Sa Earp and E. Toubiana. 
{\em Screw motion surfaces in $\hd \times \R$ and 
$\sd \times \R$}, Illinois J. Math. \textbf{49} (2005), 1323--1362.


\bibitem{[Isabel]} I. Fernandez and P. Mira. 
{\em Harmonic maps and constant mean curvature surfaces
in $\hd \times \R$}, to appear in Amer. J. of Mathematics.


\bibitem{[Pedrosa]} C. Figueroa, F. Mercuri and R. Pedrosa. 
{\em Invariant surfaces of the Heisenberg groups}, 
Ann Mat. Pura Appl. (IV), Vol 177, 173-194, 1999.



\bibitem{[Hauswirth]}  L. Hauswirth, {\em Minimal surfaces of Riemann type in
three-dimensional product manifolds}, Pacific J. Math. 224 (2006), 
91--117.

\bibitem{[Hauswirth1]} L. Hauswirth, R. Sa Earp and E. Toubiana. 
{\em Associate and conjugate minimal immersions in $\M\times \R$}, 
to appear in Tohoku Math. J.

\bibitem{[M-R]} W. Meeks III and H. Rosenberg.
{\em The theory of minimal surfaces in $M\times R$}, 
Comment. Math. Helv. (80), 811-858, 2005.





\bibitem{[Nelli-R]} B. Nelli and H. Rosenberg.
{\em Minimal surfaces in $H^2\times R$},
Bull. Braz. Math. Soc. \textbf{33}, 263-292, 2002.

\bibitem{[Rosenberg]} H. Rosenberg. 
{\em Minimal surfaces in $\M^2\times \R$}, Illinois J. Math.
 {\bf 46}, N$^o$ 4, 1177-1195, 2002. 


\bibitem{[Sanini]} A. Sanini. {\em Gauss map of a surface of the
    Heisenberg group}, Boll. Un. Mat. Ital. B (7), 11, 79-93, 1997.

\bibitem{[San-Ger]} G. Sansone and J. Gerretsen. {\em Lectures on the theory of
functions of a complex variable}, P. Noordhoff, Groningen, 1960.

\bibitem{[Scott]} P. Scott. {\em The geometries of
    3-manifolds}, Bull. London Math. Soc., 15 (5): 401-487, 1983.


\bibitem{[S-T]} R. Souam and E. Toubiana. 
{\em On the classification and regularity of umbilic surfaces in homogeneous
  3-manifolds}, Mat. Contemp., Vol 30, 201-215, 2006.




\bibitem{[Spivak]} M. Spivak. {\em A comprehensive introduction to
    differential geometry}, Vol 3 and 4, Boston, Publish or Perish, 1970.


\bibitem{[Thurston]} W. Thurston. {\em Three-Dimensional Geometry and
    Topology}, Princeton, 1997.



\end{thebibliography}
\end{document}